\documentclass[twoside,twocolumn,english,9pt,twoside,intoc,bibliography=totoc,BCOR10mm,captions=tableheading,fleqn]{article}
\usepackage[T1]{fontenc}
\usepackage[utf8]{inputenc}
\usepackage[a4paper]{geometry}
\usepackage{fancyhdr}
\usepackage{babel}
\usepackage{prettyref}
\usepackage{units}
\usepackage{algorithm2e}
\usepackage{amsmath}
\usepackage{amsthm}
\usepackage{amssymb}
\usepackage{graphicx}
\usepackage{wasysym}
\usepackage{esint}
\usepackage{nomencl}

\providecommand{\makenomenclature}{\makeglossary}
\makenomenclature
\usepackage[unicode=true,
 bookmarks=true,bookmarksnumbered=true,bookmarksopen=true,bookmarksopenlevel=1,
 breaklinks=true,pdfborder={0 0 0},backref=false,colorlinks=false]
 {hyperref}
\hypersetup{
 pdfauthor={Li, Wenwei},
 pdfpagelayout=OneColumn,pdfnewwindow=true,pdfstartview=XYZ,plainpages=false}

\makeatletter
\newcommand{\lyxaddress}[1]{
\par {\raggedright #1
\vspace{1.4em}
\noindent\par}
}

\usepackage{flushend,cuted}    

\usepackage{exscale}	
\usepackage{relsize}
\usepackage{mathcomp}  

\AtBeginDocument{
  
}

\makeatother

\begin{document}

\title{\textsf{Approximation of the Partition Number After Hardy and Ramanujan:}\\
\textsf{An Application of Data Fitting Method in Combinatorics}}

\maketitle
\begin{strip}
\begin{center}

\author{\emph{LI Wenwei}}

\lyxaddress{\noindent \begin{center}
School of Mathematical Science, University of Science and Technology
of China, \\
NO. 96, Road Jinzhai, Hefei, Anhui, P. R. China, 230026\\
E-Mail: \emph{liwenwei@ustc.edu} 
\par\end{center}}

\end{center}
\begin{abstract}
  Sometimes we need the approximate value of the partition number
in a simple and efficient way. There are already several  formulae
to calculate the partition number $p(n)$.  But they are either
inconvenient for  most people (not majored in math) who do not want
do write programs, or unsatisfying in accuracy.  By  bringing in
two parameters in the Hardy-Ramanujan's Asymptotic formula and fitting
the data of the two parameters by least square method, iteration method
and some other special designed methods, several revised elementary
estimation formulae with high accuracy for $p(n)$ are obtained. With
these estimation formulae,  the approximate value of $p(n)$ can
be calculated by a pocket calculator without programming function.
The main difficulty is that the usual methods to fit the data of the
two parameters  by an elementary function is defective here. These
method could be used in finding the fitting functions of some other
complex data.

\vspace{0.5cm}


\textbf{AMS2010 Subject Classification:} 65D10, 65D15, 05A17, 11P81,
 68-04, 68R05.
\end{abstract}
$\ $

$\ $

$\ $

\end{strip}
\maketitle\newpage{}

\flushend 

\tableofcontents{}

\newpage{}


\section{Introduction\label{sec:Introduction}}

The partition number $p(n)$ is an interesting topic which attracts
many attention. There are already a lot of literatures on many aspects
 of $p(n)$.  Many mathematicians, such as Euler, Hardy, Ramanujan,
Rademacher, Newman, Erdős, Andrews, Berndt and Ono, have made important
contribution to this topic. Some important literatures may be found
in \cite{Anonymous2007Bib-Partition}, or in the references of \cite{Eric1999PtFuncP},
\cite{KenOno2013RecentWorkPttFunc}, \cite{Stump2015IntPtt} and \cite{SloaneOEIS-A000041-pn}.

In recent years, a very important result dues to Ken Ono and his team
who connected the partition function with the modular form and found
the principles of the congruence property of $p(n)$ that may even
be considered as the revealing of the nature of numbers (refer \cite{KenOno2013News-FractalStructurePttFunc},
\cite{Amanda2012-l-adic-Properties-Elsevier}, \cite{Frank2011Remark}
 and \cite{Carol2011-News-NMTRNN}).

For a positive integer $n$, an integer solution of the equation 
\begin{equation}
s_{1}+s_{2}+\cdots+s_{q}=n\ \ (1\leqslant s_{1}\leqslant s_{2}\leqslant\cdots\leqslant s_{q},\,q\geqslant1),\label{eq:Pt-n-1}
\end{equation}
 for all the possible integer $q$ (where $s_{1}$, $s_{2}$, $\cdots$,
$s_{q}$ are unknowns) is called a \emph{partition} of  $n$. The
number of all the partitions of $n$ is denoted by $p(n)$, which
is also called the \emph{partition number} or the \emph{partition
function}.

In a lot of occasions, we need the value of $p(n)$. There are already
several formulae to calculate $p(n)$. 

In reference \cite{Marshall1958Survey} (p.53, p.57) or \cite{Ono2000ArithPtFunc},
we may find the generation function of $p(n)$ obtained by Euler:
\phantom{}
\begin{align}
F(x) & =\sum\limits _{n=0}^{\infty}p(n)x^{n}=\dfrac{1}{1-x}\dfrac{1}{1-x^{2}}\dfrac{1}{1-x^{3}}\cdots\nonumber \\
 & \dfrac{1}{1-x^{i}}\cdots\cdots=\prod\limits _{i=1}^{\infty}\left(1-x^{i}\right)^{-1},\label{eq:pn-gen-function}
\end{align}
and a formula 
\begin{equation}
p(n)=\dfrac{1}{2\pi i}\mathlarger{\oint}_{C}\dfrac{F(x)}{x^{n+1}}\mbox{d}x,\label{eq:p(n)-integral}
\end{equation}
 where $C$ is a contour around the original point. Of course, we
seldom use  \eqref{eq:p(n)-integral} to compute the value of $p(n)$
in practical.

There is a recursion for $p(n)$ (\cite{Marshall1958Survey}, p.55),
\begin{align}
p(n) & =p(n-1)+p(n-2)-p(n-5)\nonumber \\
 & \quad\quad-p(n-7)+\cdots\nonumber \\
 & \quad\quad+(-1)^{k-1}p\left(n-\dfrac{3k^{2}\pm k}{2}\right)+\cdots\cdots\nonumber \\
 & =\sum\limits _{k=1}^{k_{1}}(-1)^{k-1}p\left(n-\dfrac{3k^{2}+k}{2}\right)\nonumber \\
 & \quad\quad+\sum\limits _{k=1}^{k_{2}}(-1)^{k-1}p\left(n-\dfrac{3k^{2}-k}{2}\right),\label{eq:pn-recursion}
\end{align}
 where 
\begin{equation}
k_{1}=\left\lfloor \dfrac{\sqrt{24n+1}-1}{6}\right\rfloor ,\ k_{2}=\left\lfloor \dfrac{\sqrt{24n+1}+1}{6}\right\rfloor ,\label{eq:pn-recursion-k1-k2}
\end{equation}
 and assume that $p(0)=1$. Here $\left\lfloor x\right\rfloor $ stands
for the maximum integer that will not exceed the real number $x$.

 Equation \eqref{eq:pn-recursion} is much better for computing 
$p(n)$. We can  obtain the exact value of $p(n)$ efficiently with
a program based on it. But it is not convenient for  many people
who do not want to write programs. Further more, if we want to calculate
$p(n)$ by  \eqref{eq:pn-recursion} by a small program written in
C or some other general computer language, it is usually necessary
to decide the size of the space in memory to store the results beforehand,
which means we should know the approximate value of $p(n)$ before
the calculation started, (actually, here it is sufficient to know
$\left\lceil \frac{\log_{2}p(n)+1}{8}\right\rceil $, where $\left\lceil x\right\rceil $
stands for the minimum integer that is greater than or equal to the
real number $x$)  otherwise we have to do some extra work for overflow
handling and consequently change the size of the space in memory to
store the value of the variable that stands for $p(n)$.  

Obviously, the datatypes already defined in the C language itself
are not suitable. 

If we use the Dynamic Memory Allocation method, this problem is solved
at the price of the program being a little more complicated. Actually,
in a lot of cases, we can not decide the approximate size of the result,
it  is the best  choice available.

If we can use maple, maximal, axiom or some other computer algebra
systems, there is no need to consider this problem. But it is not
always an option, especially when the function to do this job is part
of a big program written in a compile language while mixing programming
of an interpretative language and a compile language is nearly unavailable
in most cases (with very few exceptions, such as mixing programming
C and matlab). 

$\ $

The analysis of $p(n)$ by contour integral with  \eqref{eq:p(n)-integral}
(refer \cite{Marshall1958Survey}, p. 57) resulted a very good estimation
of $p(n)$, 
\begin{equation}
p(n)=\sum\limits _{q=1}^{\left\lfloor \alpha\sqrt{n}\right\rfloor }A_{q}(n)\cdot\phi_{q}(n)+O(n^{-\nicefrac{1}{2}}),\label{eq:Hardy-Ram-fml1}
\end{equation}
 called the \emph{Hardy-Ramanujan formula \index{Hardy-Ramanujan formula}}
(refer \cite{Ramanujan1918AsymFmlComAnal} and \cite{Rademacher1937ConvergSeries}),
that 6 terms of this formula contain an error of 0.004 when $n$ =
100, while 8 terms of this formula contain an error of 0.004 when
$n$ = 200. Here $\alpha$ is an arbitrary constant, 
\[
\phi_{q}(n)=\dfrac{\sqrt{q}}{2\pi\sqrt{2}}\cdot\dfrac{\mathrm{d}}{\mathrm{d}n}\left(\dfrac{\exp\left(\tfrac{\pi}{q}\sqrt{\tfrac{2}{3}\left(n-\tfrac{1}{24}\right)}\right)}{\sqrt{n-\tfrac{1}{24}}}\right),
\]
\[
A_{q}(n)=\underset{\tiny\begin{array}{c}
0<p<q\\
(p,q)=1
\end{array}}{\mathlarger{\Sigma}}\omega_{p,q}\cdot\exp\left(\dfrac{-2np\pi i}{q}\right)
\]
 (while $p$ runs through the non-negative integers that are prime
to $q$ and less than $q$), $\omega_{p,q}$  is a certain 24$q$-th
root of unity,   $\left(\dfrac{a}{b}\right)$ is the Legendre symbol.
$b$ is an odd prime,  and $p'$ is any positive integer such that
$q\,|\,(1+pp')$. When $n$ is very large, $p(n)$ is the integer
nearest to $\sum\limits _{q=1}^{\left\lfloor \alpha\sqrt{n}\right\rfloor }A_{q}(n)\cdot\phi_{q}(n)$.

In \cite{Marshall1958Survey} or \cite{Rademacher1937ConvergSeries},
 a convergent series for $p(n)$ modified from  \eqref{eq:Hardy-Ram-fml1}
by Rademacher in 1937 is presented, 
\begin{equation}
p(n)=\sum\limits _{q=1}^{\infty}A_{q}(n)\cdot\psi_{q}(n),\label{eq:Rademacher-1}
\end{equation}
 where $A_{q}(n)$ is the same as mentioned above and
\[
\psi_{q}(n)=\dfrac{\sqrt{q}}{\pi\sqrt{2}}\dfrac{\mathrm{d}}{\mathrm{d}n}\left(\dfrac{\sinh\left(\tfrac{\pi}{q}\sqrt{\frac{2}{3}\left(n-\tfrac{1}{24}\right)}\right)}{\sqrt{n-\tfrac{1}{24}}}\right).
\]

\begin{figure}
\noindent \begin{centering}
\includegraphics[scale=0.5]{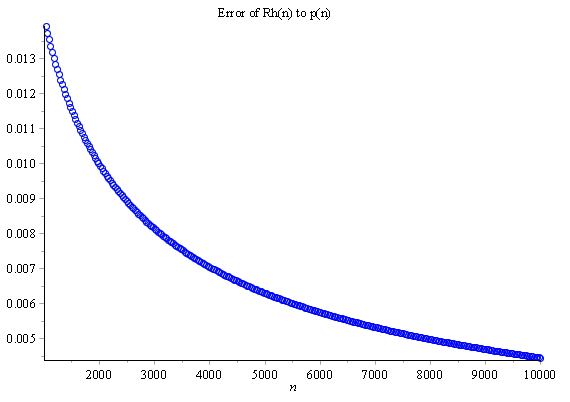}
\par\end{centering}

\noindent \begin{centering}
\label{Fig:2_4_1_0_Rel_Err_Rh_10K}
\par\end{centering}

\noindent \centering{}\caption{The Relative Error of $R_{\mathrm{h}}(n)$ to $p(n)$ when 1K $\leqslant n\leqslant$
10K}
\end{figure}

Equations \eqref{eq:Hardy-Ram-fml1} or  \eqref{eq:Rademacher-1}
are valuable in theory and can be used to calculate the value of $p(n)$
with very high accuracy. But they are not convenient for practical
usage especially when $n$ is small, since it is very difficult for
programmers, engineers or other ordinary people (not familiar with
any computer algebra system softwares) since they are too complicated
and they contain some special functions that most people (not majored
in mathematics)  are not familiar with.  It is very difficult for
them to use these two formulae to calculate $p(n)$ on a pocket science
calculator without programming function.

In references \cite{Eric1999PtFuncP} or \cite{NIST2015FNTANTUP},
we may find the famous asymptotic formula for $p(n)$, \label{Sym:Rh(n)}
\nomenclature[Rh(n)]{$R_{\mathrm{h}}(n)$}{The Hardy-Ramanujan's asymptotic formula. \pageref{Sym:Rh(n)}}
\begin{equation}
p(n)\sim\dfrac{1}{4n\sqrt{3}}\exp\left(\sqrt{\frac{2}{3}}\pi n^{\nicefrac{1}{2}}\right),\label{eq:Ram-Hardy-Fml-Est-1}
\end{equation}
 obtained by Godfrey Harold Hardy and Srinivasa Ramanujan in 1918
in the famous paper \cite{Ramanujan1918AsymFmlComAnal}. (Two different
proofs can be found in \cite{Pal1942ElemProofRHFml} and \cite{Newman1962SimProofPtFml}.
The evaluation of the constants was shown in \cite{Donald1951EvalfConstHRFml}.)
This formula will  be called the \emph{Hardy-Ramanujan's asymptotic
formula\index{Hardy-Ramanujan asymptotic formula}} in this paper.
This asymptotic formula is with great importance in theory.  Equation
\eqref{eq:Ram-Hardy-Fml-Est-1} is much more convenient than formulae
 \eqref{eq:Hardy-Ram-fml1} and  \eqref{eq:Rademacher-1} for ordinary
people  not majored in mathematics. 

Let 
\begin{equation}
R_{\mathrm{h}}(n)=\dfrac{1}{4n\sqrt{3}}\exp\left(\sqrt{\frac{2}{3}}\pi\sqrt{n}\right).\label{eq:Ram-Hardy-Formula}
\end{equation}

be the asymptotic function by Hardy and Ramanujan. 

By the figure in reference \cite{SloaneOEIS-A002865graph}, this asymptotic
formula fits $p(n)$ very well when $n$ is huge. But when $n$ is
small, the relative error of $R_{\mathrm{h}}(n)$ to $p(n)$ is not
so satisfying as shown in Table \ref{Table:Rel-Err-p(n)-HR(n)} (when
$n$ $\leqslant$ 1000) on page \pageref{Table:Rel-Err-p(n)-HR(n)}.
When $n$ $\leqslant$ 25, the relative error is greater than 9\%;
when 25 < $n$ $\leqslant$ 220, the relative error is greater than
3\%;  when $n$ $\leqslant$ 1000, the relative error is greater
than 1.4\%. From Figure 1.\ref{Fig:2_4_1_0_Rel_Err_Rh_10K}, we will
find out that the relative error is greater than 0.44\% when 1000
$\leqslant n\leqslant$ 10000. Considering that $p(n)$ is an integer
and $R_{\mathrm{h}}(n)$ is definitely not, the round approximation
of $R_{\mathrm{h}}(n)$ may be a little more accurate, but that does
not help.

Although  \eqref{eq:Hardy-Ram-fml1} is not so accurate when $n$
is small, it provides some important clue for a  more accurate formula
for small $n$. 

\begin{table}
\noindent 

\noindent \begin{centering}
\includegraphics[bb=97bp 409bp 435bp 597bp,clip,scale=0.63]{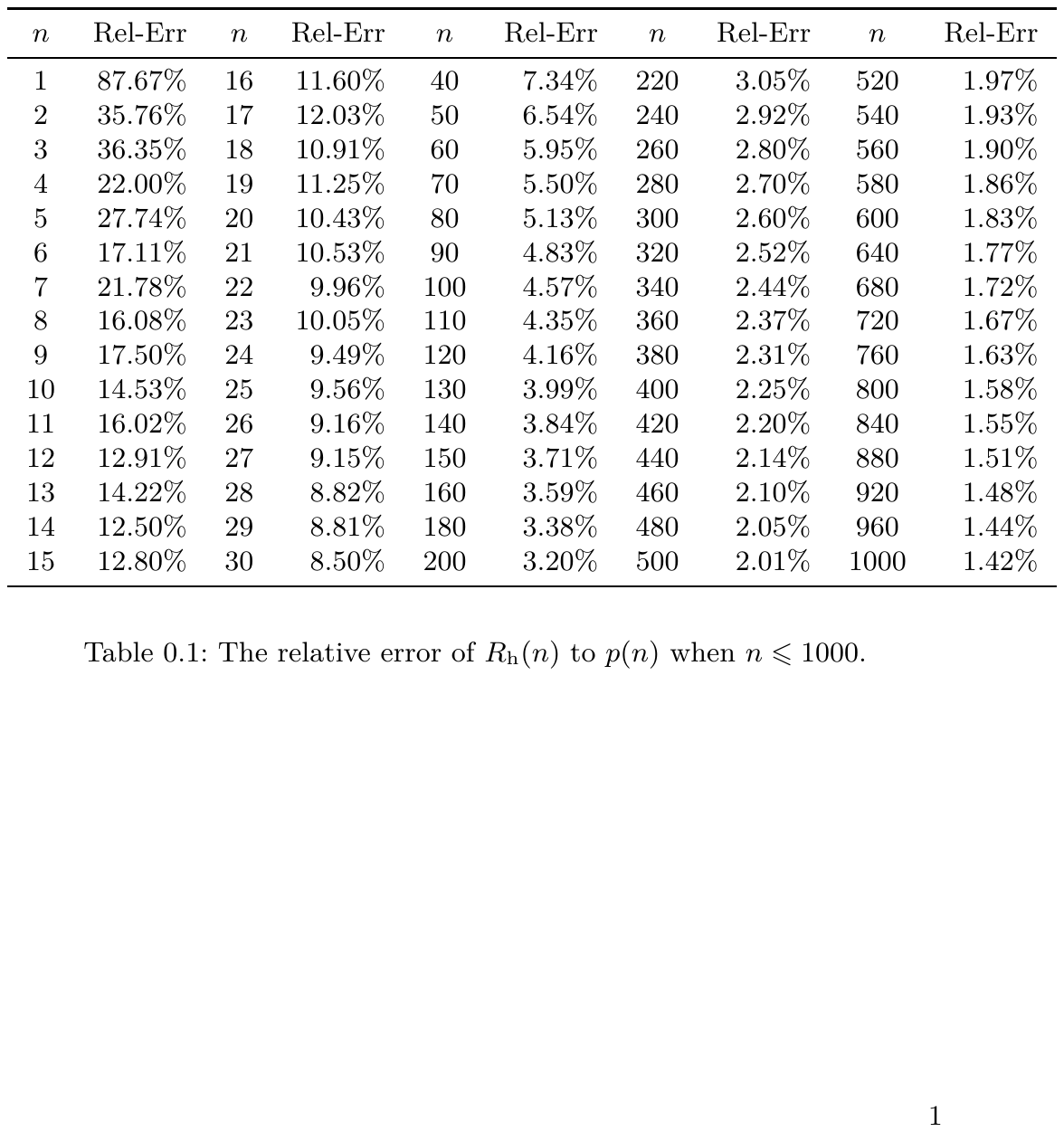}
\par\end{centering}

\caption{The relative error of $R_{\mathrm{h}}(n)$ to $p(n)$ when $n\leqslant1000$.}
\label{Table:Rel-Err-p(n)-HR(n)}

\end{table}

By revising \eqref{eq:Hardy-Ram-fml1}, some other estimation formulae
with high accuracy is obtained here. 

In \prettyref{sec:main-idea}, the main idea is introduced, two parameters
$C_{1}$ and $C_{2}$ are brought in the Hardy-Ramanujan's asymptotic
formula, they will be fitted in sections \ref{sec:Fit-the-Exponent-C1(n)}
and \ref{sec:Fit-the-Denominator-C2(n)}, respectively. Sections \ref{sec:Estmt-p(n)-other-methd}
and \ref{sec:Estimate-p(n)-by-fitting-Rh(n)-p(n)} will show some
other methods to obtain estimation formulae. Section \ref{sec:Estmt-p(n)-less-than-100}
displays an estimation formula with more accuracy when $n\leqslant100$. 

The main difficulty is that it is too hard to obtain the appropriate
functions to fit the data of $C_{1}$ (or $C_{2}$ or some others)
generated here since we know very little about them and the usual
methods to find fitting functions are invalid here. If we fit the
data directly, the results are far from satisfactory, at least the
accuracy is not as good as that of \eqref{eq:Ram-Hardy-Fml-Est-1}.

\section{Main idea \label{sec:main-idea}}

There are many different ways to modify $R_{\mathrm{h}}(n)$, e.g.
we could also construct a function $p_{1}(n)$ to estimate $R_{\mathrm{h}}(n)-p(n)$,
then $R_{\mathrm{h}}(n)-p_{1}(n)$ may reach a better accuracy when
estimating $p(n)$, or we can estimate the value of $\frac{R_{\mathrm{h}}(n)}{p(n)}$
by a function $f_{1}(n)$ then estimate $p(n)$ by $\frac{R_{\mathrm{h}}(n)}{f_{1}(n)}$,
etc. The problem is that the accuracy of $R_{\mathrm{h}}(n)-p_{1}(n)$
is not so satisfying if we do not use the idea shown in  \eqref{eq:Ram-Hardy-Rev-Form1}
in \prettyref{sec:main-idea},  because the shape of the figure of
$\ln\left(R_{\mathrm{h}}(n)-p(n)\right)$ is nearly the same as the
shape of the figure of $\ln\left(p(n)\right)$, at least we can not
tell the difference  by our eyes as shown on Figure \ref{Fig:2_4_19_ln(p(n))}
and Figure \ref{Fig:2_4_20_ln(Rh(n)-p(n))} (on page \pageref{Fig:2_4_20_ln(Rh(n)-p(n))}),
though they are different in theory.

$\ $

Since $p(n)$ $\sim$ $R_{\mathrm{h}}(n)$,  we believe that an approximate
formula with better accuracy may be in this form 
\begin{equation}
p(n)\approx\dfrac{1}{4\sqrt{3}(n+C_{2})}\exp\left(\sqrt{\frac{2}{3}}\pi\sqrt{n+C_{1}}\right).\label{eq:Ram-Hardy-Rev-Form1}
\end{equation}

\begin{figure}[h]
\centering\includegraphics[scale=0.29]{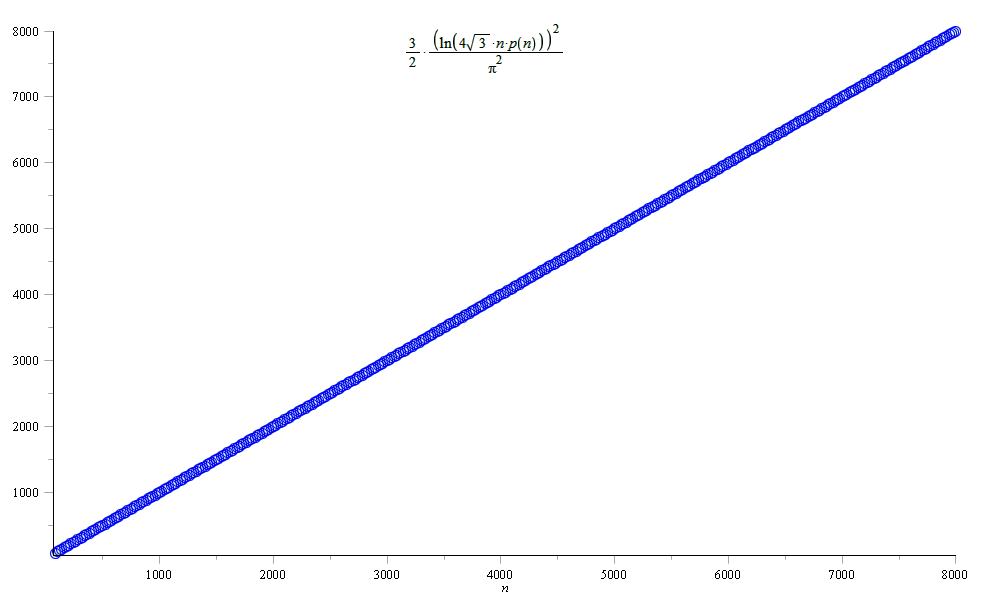}\caption{The graph of the data  $\left(n,\ \frac{3}{2}\cdot\frac{\left(\ln\left(4n\sqrt{3}p(n)\right)\right)^{2}}{\pi^{2}}\right)$.}
\label{Fig:2_4_1_(n_C1(n)+n)}

\medskip{}

\centering\includegraphics[scale=0.28]{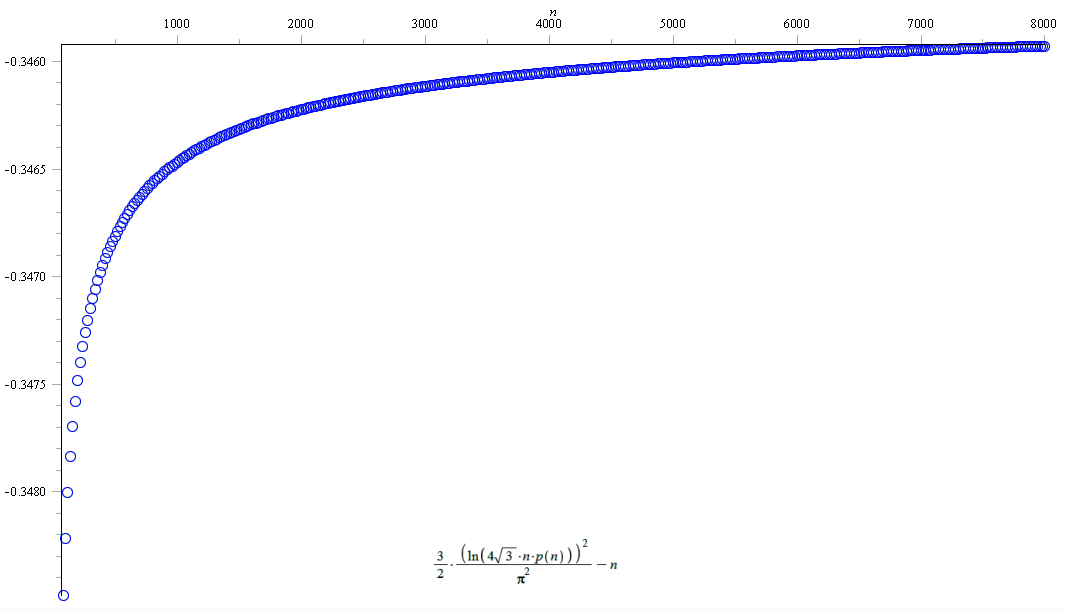}\caption{The graph of the data $\left(n,\,C_{1}\left(n\right)\right)$ ($n\geqslant120$).}
\label{Fig:2_4_2_(n_C1(n))}
\end{figure}

Where $C_{1}$ (or $C_{2}$) may be a constant or a function of $n$
that increases slowly than $n$, so as to have $\underset{n\rightarrow\infty}{\lim}\frac{\frac{1}{4\sqrt{3}(n+C_{2})}\exp\left(\sqrt{\frac{2}{3}}\pi\sqrt{n+C_{1}}\right)}{R_{\mathrm{h}}(n)}$
= 1, or $\underset{n\rightarrow\infty}{\lim}\frac{\frac{1}{4\sqrt{3}(n+C_{2})}\exp\left(\sqrt{\frac{2}{3}}\pi\sqrt{n+C_{1}}\right)}{p(n)}$
= 1.

There are some other ways to modify $R_{\mathrm{h}}(n)$, we will
discuss the details in section \ref{sec:Estmt-p(n)-other-methd}.

As we can not determine $C_{1}$ and $C_{2}$ at the same time because
of technique problems, \footnote{$\ $ Usually, we will get the value of $C_{1}$ and/or $C_{2}$ from
a number of pairs of $\left(n,p(n)\right)$ by the least square method,
not from two pairs of $\left(n,p(n)\right)$ only. Many software can
get efficiently the undetermined coefficients (by the least square
method) by solving a system of (incompatible) linear equations, while
it is very difficult to ``solve'' a system of tens or hundreds of
transcendental equations that are incompatible.} we may decide $C_{1}$ first then determine $C_{2}$, the main reason
is that $\frac{1}{(n+C_{2})}$ and $\frac{1}{n}$ differs very little
when $n$ is very huge, at least we believe that the difference is
much less that the difference of $\exp\left(\sqrt{\frac{2}{3}}\pi\sqrt{n+C_{1}}\right)$
and $\exp\left(\sqrt{\frac{2}{3}}\pi\sqrt{n}\right)$. \footnote{$\ $ It is not difficult to know that $\frac{1}{(n+\delta)}$ $\approx$
$\frac{1}{n}\left(1-\frac{\delta}{n}\right)$, $\exp\left(\sqrt{\frac{2}{3}}\pi\sqrt{n+\delta}\right)$
$\approx$ $\exp\left(\sqrt{\frac{2}{3}}\pi\sqrt{n}\right)\left(1+\frac{\pi}{\sqrt{6}}\frac{\delta}{\sqrt{n}}\right)$,
when $\delta\ll n$. Obviously, $\frac{\delta}{n}\ll\frac{\pi}{\sqrt{6}}\frac{\delta}{\sqrt{n}}$
(when $\max\{\delta,1\}\ll n$).} 

\begin{figure}[h]
\centering\includegraphics[scale=0.31]{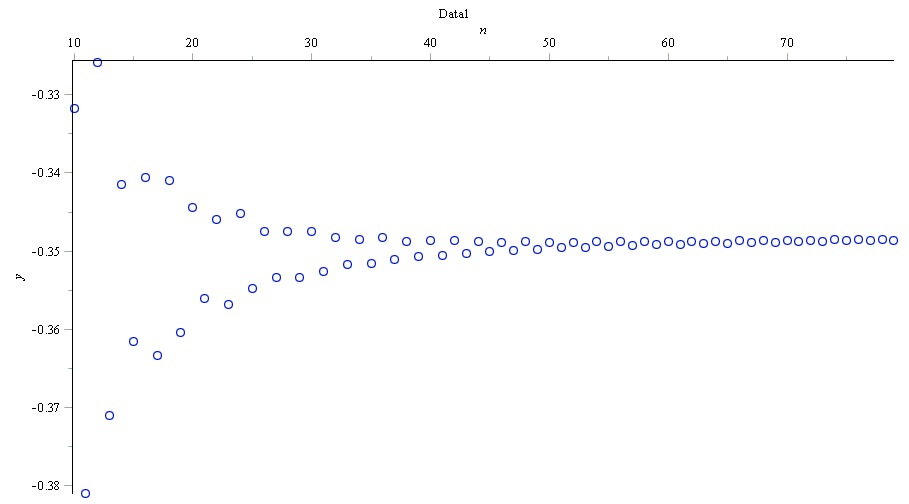}\caption{The graph of the data $\left(n,\ \frac{3}{2}\cdot\frac{\left(\ln\left(4n\sqrt{3}p(n)\right)\right)^{2}}{\pi^{2}}\right)$
($n\leqslant80$).}
\label{Fig:2_4_3_(n_C1(n))-part-1}

$\ $

\centering\includegraphics[scale=0.31]{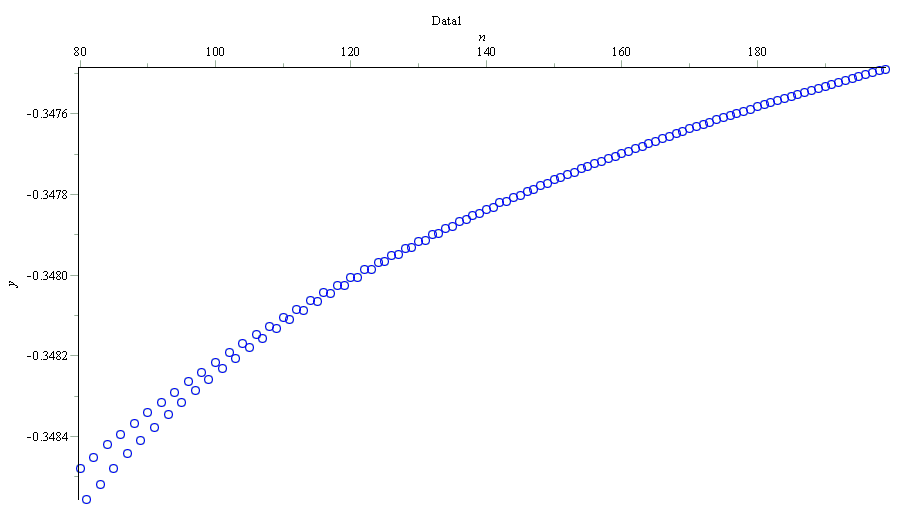}\caption{The graph of the data $\left(n,\,C_{1}\left(n\right)\right)$ ($80\leqslant n\leqslant200$).}
\label{Fig:2_4_4_(n_C1(n))-part-2}
\end{figure}

So, when $n$ $\gg$ 1, we believe 
\[
p(n)\doteq\dfrac{1}{4\sqrt{3}n}\exp\left(\sqrt{\frac{2}{3}}\pi\sqrt{n+C_{1}}\right),
\]

hence $4\sqrt{3}n\times p(n)\doteq\exp\left(\pi\sqrt{\frac{2}{3}(n+C_{1})}\right)$,
then 
\begin{equation}
C_{1}\left(n\right)\doteq\dfrac{3}{2}\cdot\dfrac{\left(\ln\left(4n\sqrt{3}p(n)\right)\right)^{2}}{\pi^{2}}-n.\label{eq:C1-deduce-fml-approximate}
\end{equation}
If we point the data $\left(n,\ \frac{3}{2}\cdot\frac{\left(\ln\left(4n\sqrt{3}p(n)\right)\right)^{2}}{\pi^{2}}\right)$
($n=20k+100$, $k$ = 1, 2, $\cdots$, 395) in the coordinate system,
we will find that they lie in a straight line, as shown in the  Figure
\ref{Fig:2_4_1_(n_C1(n)+n)} on page \pageref{Fig:2_4_1_(n_C1(n)+n)},
which means that the Hardy-Ramanujan's asymptotic formula is close
to perfect. Here every tiny cycle stands for a data point.

\section{Fit the Exponent \label{sec:Fit-the-Exponent-C1(n)}}

\begin{figure}[h]

\noindent \begin{centering}
\includegraphics[scale=0.26]{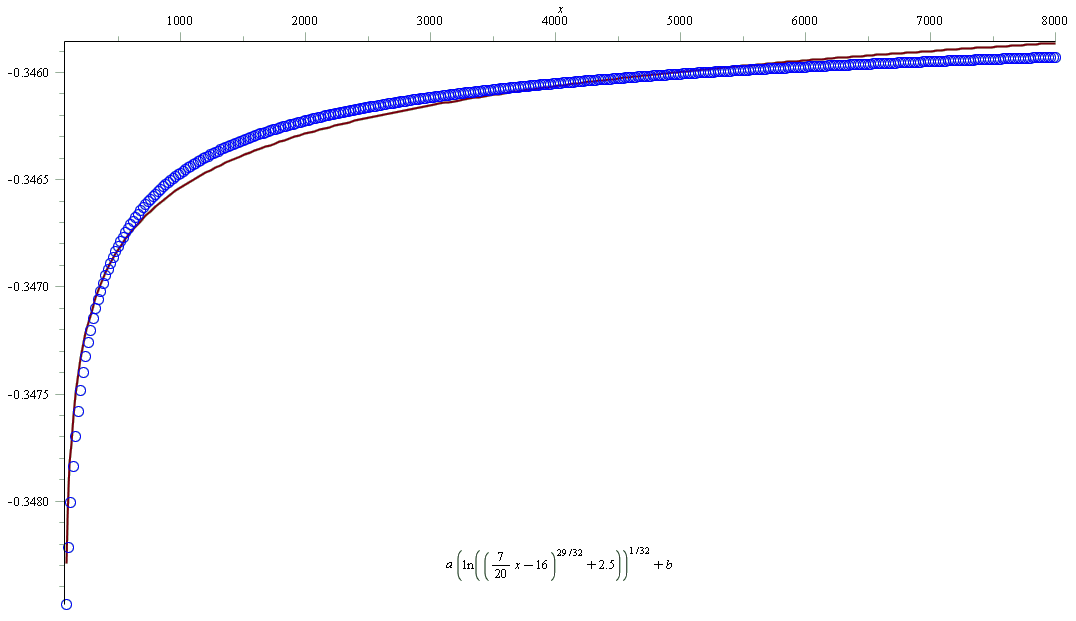}
\par\end{centering}

\noindent \begin{centering}
\caption{The graph of a bad fitting curve of the data $\left(n,\,C_{1}\left(n\right)\right)$}
\label{Fig:2_4_5_(n_h(n))-first-part-bad}
\par\end{centering}

\end{figure}

If we point the data $\left(n,\,C_{1}(n)\right)$, i.e., $\left(n,\ \frac{3}{2}\cdot\frac{\left(\ln\left(4n\sqrt{3}p(n)\right)\right)^{2}}{\pi^{2}}-n\right)$
($n=20k+100$, $k$ = 1, 2, $\cdots$, 395) in the coordinate system,
we will get the  Figure \ref{Fig:2_4_2_(n_C1(n))} on page \pageref{Fig:2_4_2_(n_C1(n))}.
Here the points when $n\leqslant120$ are not shown on  Figure \ref{Fig:2_4_2_(n_C1(n))},
partly because the deduction above is based on $n\gg1$, the main
reason is that the points obviously do not lie in a curve when $n\leqslant120$,
as shown on Figure \ref{Fig:2_4_3_(n_C1(n))-part-1} and Figure \ref{Fig:2_4_4_(n_C1(n))-part-2}
(on page \pageref{Fig:2_4_4_(n_C1(n))-part-2}).

 Figure \ref{Fig:2_4_2_(n_C1(n))} looks like a logarithmic curve
or a hyperbola. The author has tried many functions (by a small program
written in MAPLE) like
\[
a\cdot\left(\ln(x^{e_{1}}+c_{1})\right)^{e_{2}}+b,
\]
 where $e_{1}$, $e_{2}$ and $c_{1}$ are given constants while $a$
and $b$ are undermined coefficients to be decided.  But none of
them fits the data very well. A function 
\[
y=a\cdot\left(\ln\left(\left(\dfrac{7}{20}\cdot x-16\right)^{29/32}+2.5\right)\right)^{1/32}+b,
\]
 where $a$ = 0.06656839293 and $b$ = -0.4166945066, may fit the
data better, but it is not as good as we expect, as shown on  Figure
\ref{Fig:2_4_5_(n_h(n))-first-part-bad} on page \pageref{Fig:2_4_5_(n_h(n))-first-part-bad}.

A hyperbola like $y=\dfrac{a}{x}+b$ does not fit the data very well,
either. Then we consider this type of functions 
\begin{equation}
y=\dfrac{a}{(x+c_{2})^{e_{2}}}+b,\label{eq:C1-Fit-type}
\end{equation}
 where $a$, $b$, $c_{2}$ and $e_{2}$ are undetermined constants.
This seems much better. For technique reason, we can not decide all
the undetermined coefficients $a$, $b$, $c_{2}$, $e_{2}$ at the
same time. \footnote{$\ $ Because most computer algebra system (CAS) could not solve the
system of many incompatible nonlinear equations by the least square
method, or the time-consumption is unacceptable. }

\begin{figure}[h]
\noindent \begin{centering}
\centering\includegraphics[scale=0.26]{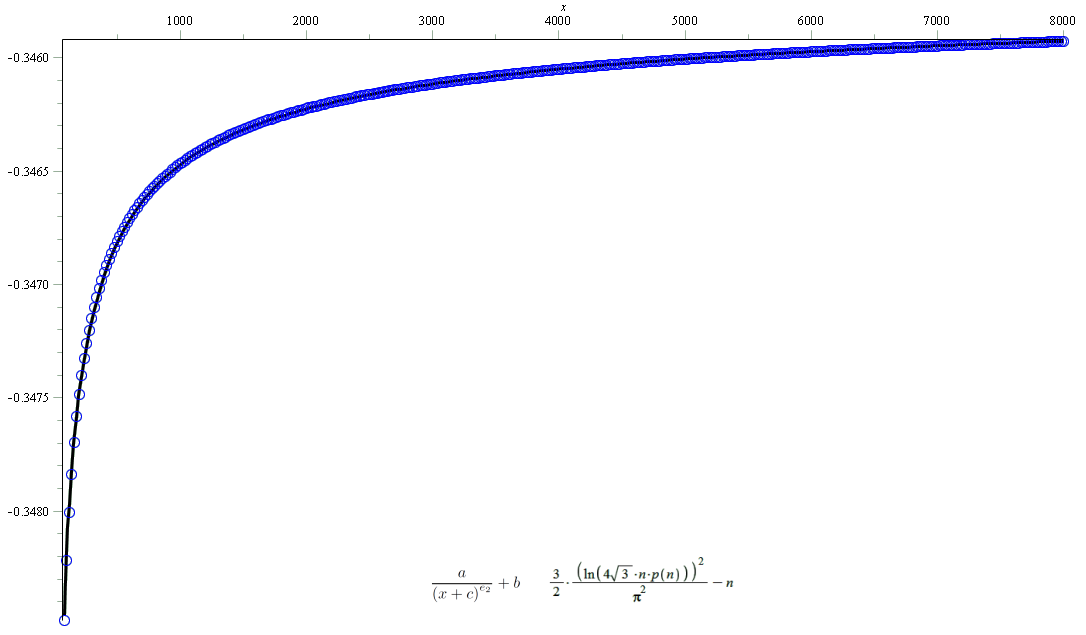}
\par\end{centering}

\noindent \begin{centering}
\caption{The graph of a good fitting curve of the data $\left(n,\,C_{1}\left(n\right)\right)$}
\label{Fig:2_4_6_(n_C1(n))-fit_1}
\par\end{centering}

\noindent \begin{centering}
$\ $\medskip{}

\par\end{centering}

\noindent \begin{centering}
\centering\includegraphics[scale=0.28]{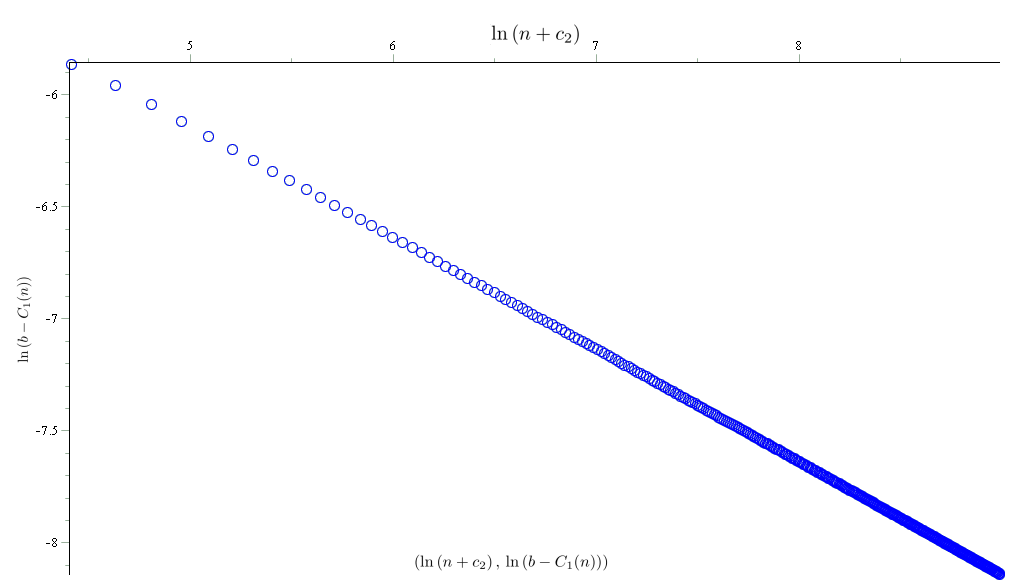}
\par\end{centering}

\noindent \centering{}\caption{The graph of the data $\bigl(\ln\left(n+c_{2}\right),\,\ln\left(b-C_{1}(n)\right)\bigr)$}
\label{Fig:2_4_7_(ln(n+c2)_ln(b-C1))}
\end{figure}

\begin{figure}[h]
\centering\includegraphics[scale=0.26]{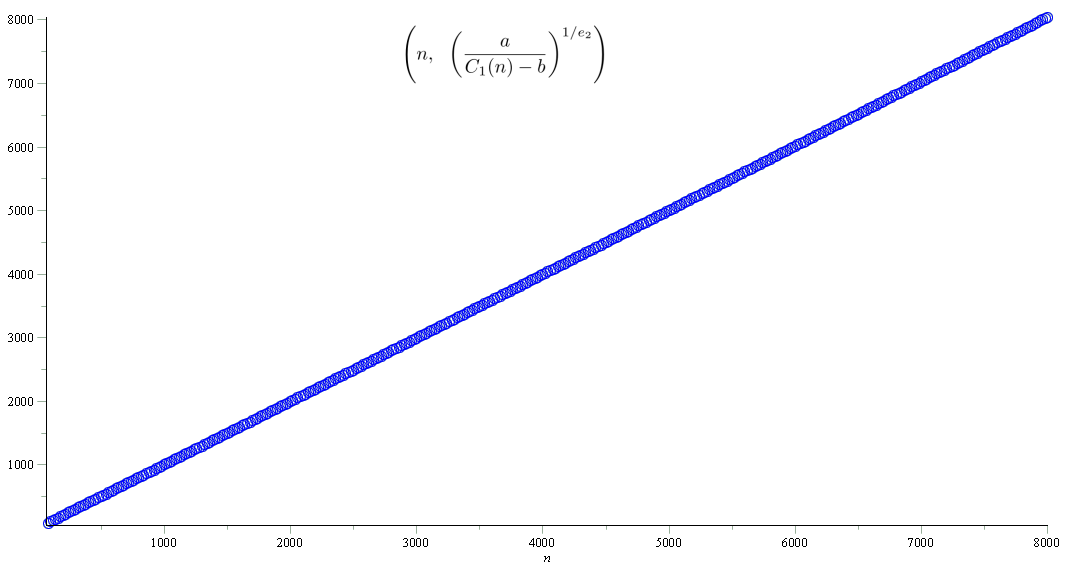}

\caption{The graph of the data $\left(n,\,\cdot\left(\frac{a}{C_{1}(n)-b}\right)^{1/e_{2}}\right)$
}
\label{Fig:2_4_8_(n_n+c2)}

$\ $

\centering\includegraphics[scale=0.26]{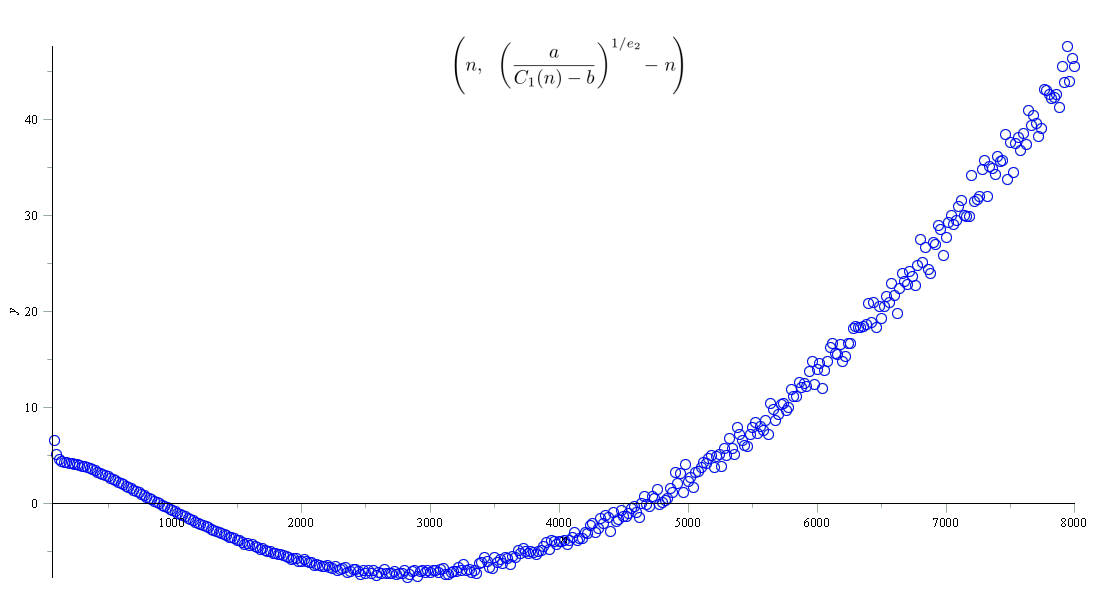}

\caption{The graph of the data $\left(n,\,\cdot\left(\frac{a}{C_{1}(n)-b}\right)^{1/e_{2}}-n\right)$
}
\label{Fig:2_4_9_(n,)}
\end{figure}

These undetermined coefficients may be obtained in this way: 
\begin{itemize}
\item A1. Give $c_{2}$ and $e_{2}$ initial values;
\item A2. Fit the data $\left(n,\,C_{1}(n)\right)$ by the least square
method with Equation  \eqref{eq:C1-Fit-type} and obtain the values
of $a$ and $b$, then get the average error of the fitting function
for the values of $c_{2}$, $e_{2}$, $a$, $b$; \footnote{$\ $ Here we use the square root of the mean square deviation 
\[
s=\sqrt{\dfrac{1}{m}\sum_{i=1}^{m}\left(y_{i}-f(x_{i})\right)^{2}}
\]
to measure the average error of the fitting function $y=f(x)$ to
the original data $\left(x_{i},\,y_{i}\right)$ ($i$ = 1, 2, $\cdots$,
$m$).  }
\item A3. Reevaluate $e_{2}$ and $a$. Plot the points of the data $\Bigl(\ln\left(n+c_{2}\right),\,\ln\left(b-C_{1}(n)\right)\Bigr)$
($n=20k+100$, $k$ = 1, 2, $\cdots$, 395) in the coordinate system
with the values of $b$ and $c_{2}$ just found, \footnote{$\ $ Such as shown in  Figure \ref{Fig:2_4_7_(ln(n+c2)_ln(b-C1))}
on page \pageref{Fig:2_4_7_(ln(n+c2)_ln(b-C1))} when $c_{2}$ = 2.5
and $b$ = $-$0.3456348045. 

$\ $ The purpose of this step is to obtain more accurate values of
$e_{2}$ and $a$. Since $C_{1}(n)=\frac{a}{(n+c_{2})^{e_{2}}}+b$,
then $b-C_{1}(n)=\frac{-a}{(n+c_{2})^{e_{2}}}$, (considering that
$a<0$), $\ln\left(b-C_{1}(n)\right)$ = $\ln\left(-a\right)$ + $e_{2}\cdot\ln\left(n+c_{2}\right)$,
so the figure of data $\Bigl(\ln\left(n+c_{2}\right),\,\ln\left(b-C_{1}(n)\right)\Bigr)$
will be some points on a straight line if the previous assumption
is correct and meanwhile the values of $b$ and $c_{2}$ are proper.} fit the data by the least square method with 
\[
y=e_{1}\cdot x+a_{1}
\]
 and find the values of $a_{1}$ and $e_{1}$, then reevaluate $e_{2}$
and $a$ by 
\[
e_{2}=-e_{1},\ a=-\exp(a_{1});
\]

\item A4. Reevaluate $c_{2}$. Plot the points of the data $\left(n,\,\left(\frac{a}{C_{1}(n)-b}\right)^{1/e_{2}}\right)$
($n=20k+100$, $k$ = 1, 2, $\cdots$, 395) in the coordinate system
with the value of $b$ and the new values of $a$ and $e_{2}$, \footnote{$\ $ Such as shown on  Figure \ref{Fig:2_4_8_(n_n+c2)} on page
\pageref{Fig:2_4_8_(n_n+c2)} when $b$ = $-$0.3456365954, $e_{2}$
= 0.5012314726 and $a$ = $-$0.02661232627.

$\ $ The main idea of this step: since $C_{1}(n)=\frac{a}{(n+c_{2})^{e_{2}}}+b$,
then $n+c_{2}=\left(\frac{a}{C_{1}(n)-b}\right)^{1/e_{2}}$, hence
the figure of data $\left(n,\,\cdot\left(\frac{a}{C_{1}(n)-b}\right)^{1/e_{2}}\right)$
will be some points on a straight line.} fit the data by the least square method with 
\[
y=x+c_{1}
\]
 and find the value of $c_{1}$, then reevaluate $c_{2}$ by $c_{2}=c_{1}$. 
\item A5. goto step 2 until a fitting function with the least  average
error is obtained. 
\end{itemize}

For example, in step A1, the initial value could be set by $c_{2}$
= 2.5, $e_{2}$ = 0.5 (or some other values).

In step A2, if $c_{2}$ = 2.5, $e_{2}$ = 0.5, then $a$ = $-$0.02635983935,
$b$ = $-$0.3456348045. \\
$\ $ If we plot the figure of  \eqref{eq:C1-Fit-type} with the
value of $c_{2}$, $e_{2}$, $a$, $b$, and compare the figure with
 Figure \ref{Fig:2_4_2_(n_C1(n))} on page \pageref{Fig:2_4_2_(n_C1(n))},
we will get a graph nearly the same as  Figure \ref{Fig:2_4_2_(n_C1(n))}
(although there should be a little different, but we can not distinguish
the difference by our eyes). The average error of the fitting function
for the values of $c_{2}$, $e_{2}$, $a$, $b$ mentioned above is
1.074574171$\times10^{-5}$, which seems to be very tiny.

In Step A3, $\ $ if $c_{2}$ = 2.5, $b$ = $-$0.3456365954, then
$a_{1}$ = $-$3.626380777, $e_{1}$ = $-$0.5012314726. \\
After reevaluation,  $e_{2}$ = 0.5012314726, $a$ = $-$0.02661232627.

In Step A4, for the values of $b$, $e_{2}$ and $a$ mentioned before,
after reevaluation $c_{2}$ = 4.871833842.

Actually, only a few times of repeating the steps form A2 to A4, we
will obtain a very good fitting function, as shown on  Figure \ref{Fig:2_4_6_(n_C1(n))-fit_1}
on page \pageref{Fig:2_4_6_(n_C1(n))-fit_1}.

$\ $  For the initial value $c_{2}$ = 2.5, $e_{2}$ = 0.5, after
repeating 41 times of the steps from A2 to A4, we will find a fitting
function 
\begin{equation}
y=\dfrac{-0.02594609078}{(x+3.320623832)^{0.4963284361}}-0.3456286995,\label{eq:C1(n)_1-A}
\end{equation}
 with a minimal average error $9.010349470\times10^{-8}$. After a
few times more of iteration, a result with similar coefficients will
be found but with a little more error.

There are some explanations about the steps above:

\begin{itemize}
\item \textsf{(1).}\textsf{\textbf{ }}In step A4, we did not plot the points
of the data $\left(n,\,\left(\frac{a}{C_{1}(n)-b}\right)^{1/e_{2}}-n\right)$
because the shape of the figure is not a horizontal line as shown
on  Figure \ref{Fig:2_4_9_(n,)} on page \pageref{Fig:2_4_9_(n,)}
(the points in the right hand side are not so smooth because only
10 significance digits are kept in the process, if more significance
digits are calculated, it will be better). Actually, it is a little
complicated. But it will not help us to obtain better values of the
undetermined in  \eqref{eq:C1-Fit-type} if we fit the data $\left(n,\,\left(\frac{a}{C_{1}(n)-b}\right)^{1/e_{2}}-n\right)$
with a more accurate fitting function. 
\item \textsf{(2). }In step A3, if we do not reevaluate $a$, the fitting
parameters will not converge in general (even if we computing more
significant figures in the process), or we can not continue the iterations
steps at all since imaginary numbers appear. 
\item \textsf{(3). }If we started with a different initial value of $c_{2}$
and keep the initial value of $e_{2}$, such as $c_{2}$ = 15, after
repeating 78 times of the steps from A2 to A4, we will find a fitting
function 
\begin{equation}
y=\dfrac{-0.02593608938}{(x+3.272445238)^{0.4962730054}}-0.3456286681,\label{eq:C1(n)_1-B}
\end{equation}
 with a minimal average error $9.109686836\times10^{-8}$.\\
\textsf{ }If we started with some different initial values for both
$c_{2}$ and $e_{2}$, such as $c_{2}$ = 15 and $e_{2}$ = 0.7, (from
 Figure \ref{Fig:2_4_2_(n_C1(n))} on page \pageref{Fig:2_4_2_(n_C1(n))},
we will find that $e_{2}$ should be less that 1.0), we will get a
similar result. After repeating 125 times of the steps from A2 to
A4, we will find a fitting function 
\begin{equation}
y=\dfrac{-0.02593617719}{(x+3.273513225)^{0.4962727258}}-0.3456286655,\label{eq:C1(n)_1-C}
\end{equation}
 with a minimal average error $9.105941452\times10^{-8}$. After that,
$e_{2}$ and $c_{2}$ will decrease slowly and slowly, and the average
error will increase little by little if we continue the steps from
A2 to A4. 
\end{itemize}
As concerned to the errors in computing, the valid value of the undermined
$a$, $b$, $c_{2}$ and $e_{2}$ should be $-0.0259361$, $-0.34562866$,
$3.273$, $0.49627$, the average absolute error of the fitting function
of $C_{1}(n)$ is about $9.1\times10^{-8}$.

Considering that  \eqref{eq:C1-deduce-fml-approximate} is an approximate
formula, we may believe that the best value of $e_{2}$ is $0.5$,
as we prefer a simple exponent. Then it will be more convenient to
obtain $a$, $b$ and $c_{2}$. 

Below $e_{2}$ is  supposed to be $1/2$, which means that the fitting
function of $C_{1}(n)$ is 
\begin{equation}
y=\dfrac{a}{\sqrt{x+c_{2}}}+b.\label{eq:C1-Fit-type-new}
\end{equation}
  When $e_{2}$ is fixed to be $1/2$, if we use the iteration method
described above but keep the value of $e_{2}$ in step A3, i.e., substitute
step A3 by 

\hphantom{If if if } A3'. Reevaluate $a$ by \footnote{$\ $ or equivalently, Plot the points of the data $\Bigl(\ln\left(n+c_{2}\right),\,\ln\left(b-C_{1}(n)\right)\Bigr)$
($n=20k+100$, $k$ = 1, 2, $\cdots$, 395) in the coordinate system
with the values of $b$, $e_{2}$ and $c_{2}$ just found, fit the
data by the least square method with $y=e_{2}\cdot x+a_{1}$ and find
the values of $a_{1}$, then reevaluate $a$ by $a=-\exp(a_{1});$}
\begin{align*}
a & =-\exp\biggl(\dfrac{1}{395}\sum\limits _{k=1}^{395}\Bigl(\ln\left(b-C_{1}(20k+100)\right)-\\
 & \quad\quad e_{2}\cdot\ln\left(20k+100+c_{2}\right)\Bigr)\biggr);
\end{align*}
 (that means we evaluate $a$ twice in every loop) the sequence of
fitting functions of $C_{1}(n)$ will diverge. But we will obtain
a converged sequence of the determinants if $n$ ranges from 120 to
6000, (i.e., consider only the data $\left(n,\,p(n)\right)$ when
$n=20k+100$, $k$ = 1, 2, $\cdots$, 295). The fitting function of
$C_{1}(n)$ obtained in this way (when $n$ ranges from 120 to 6000,
step 20) is 
\begin{equation}
y=\dfrac{-0.02650620466}{\sqrt{x+4.855479108}}-0.3456326154,
\end{equation}

with the minimal average error 2.374935895$\times10^{-7}$. \footnote{$\ $ If we use the value of $c_{1}$ already found above, such as
$c_{2}$= 3.273513225 in  \eqref{eq:C1(n)_1-C}, the fitting function
is 
\[
y=\dfrac{-0.02640970103}{\sqrt{x+3.273513225}}-0.3456340228,
\]
 with an average error $7.404647856\times10^{-7}$, which is about
3 times than that above. 

If we choose $c_{2}$= 3.320623832 in  \prettyref{eq:C1(n)_1-A},
the fitting function is 
\[
y=\dfrac{-0.02641281526}{\sqrt{x+3.320623832}}-0.3456339736,
\]
 with an average error $7.205944166\times10^{-7}$.}

For the fixed value $1/2$ of $e_{2}$, if we continue use the iteration
method described above but ignore step 3, which means we reevaluate
$a$ only once in every loop, we will meet the same situation. The
sequence of fitting functions of $C_{1}(n)$ will diverge if $n$
ranges from 120 to 8000 (or 6000) even if we calculate more significance
digits (such as 18 significance digits) in the process, but it will
converge if $n$ ranges from 120 to 4000. The fitting function of
$C_{1}(n)$ obtained in this way (when $n$ ranges from 120 to 4000,
step 20) is 
\begin{equation}
y=\dfrac{-0.02647712648}{\sqrt{x+4.55083607}}-0.345633305,
\end{equation}

with the minimal average error 1.993012726$\times10^{-7}$ when the
initial value of $c_{2}$ is 10 (iterated 4 times). But after more
times of iteration, for several initial values of $c_{2}$ (such as
5, 10, 15, etc), the fitting functions converge to 
\begin{equation}
y=\dfrac{-0.0268\cdots}{\sqrt{x+4.888\cdots}}-0.345632760\cdots,
\end{equation}

with the average error 2.68$\cdots\times10^{-7}$. $\ $

Unlike the previous method, by the results mentioned above and some
other results not mentioned here, the sequence of fitting functions
of $C_{1}(n)$ usually converges to a function which is obviously
different from the one with the minimal average error. 

In order to get a fitting function with errors as tiny as possible,
we can design another algorithm. 

By the results described above, we known that $c_{2}$ is probably
between 3 and 5, so we can find the fitting function of $C_{1}(n)$
and the corresponding average error for many values of $c_{2}$ in
the possible range, then choose the one with minimal average error.
To be cautious, we test the value of $c_{2}$ in the interval $[0.5,\ 15]$.
The main steps are as below: 
\begin{itemize}
\item (1) Initial $c_{\mathrm{a}}$, $c_{\mathrm{b}}$, $c_{0}$, $s_{0}$,
$D_{\mathrm{t}}$, $a_{0}$, $b_{0}$. Let $c_{\mathrm{a}}$ $:=$
0.5, $c_{\mathrm{b}}$ $:=$ 15, $c_{0}$ $:=$ 0, $s_{0}$ $:=$
1, $a_{0}$ $:=$ 0, $b_{0}$ $:=$ 0, $D_{\mathrm{t}}$ $:=$ 8,
$s_{\mathrm{t}}$ $:=$ 0.1,.
\item (2) for $c_{2}$ from $c_{\mathrm{a}}$ to $c_{\mathrm{b}}$ by $s_{\mathrm{t}}$
do\\
\hphantom{WW } Fit the data $\left(n,\,C_{1}(n)\right)$ by the least
square method with  \eqref{eq:C1-Fit-type-new} and get the values
of $a$ and $b$, then get the average error $s_{1}$ of the fitting
function for the values of $c_{2}$, $a$, $b$;\\
if $s_{1}$ < $s_{0}$, then let $c_{0}$ $:=$ $c_{2}$, $s_{0}$
$:=$ $s_{1}$, $a_{0}$ $:=$ $a$, $b_{0}$ $:=$ $b$; end if;
\\
end do
\item (3) If $D_{\mathrm{t}}$ > 1, then set $D_{\mathrm{t}}$ $:=$ $D_{\mathrm{t}}-1$,
$c_{\mathrm{a}}$ $:=$ $c_{0}-5s_{\mathrm{t}}$, $c_{\mathrm{b}}$
$:=$ $c_{0}+5s_{\mathrm{t}}$; \\
\hphantom{WW } set $s_{\mathrm{t}}$ $:=$ $s_{\mathrm{t}}/10$;
goto step (2);\\
else, terminate the process.\\
end if;\label{Alg:gen-c2-by-loop}
\end{itemize}
Here the symbol ``$x$ $:=$ $y$'' means that the variable $x$
is evaluated by the value of the variable $y$; in step (1), $D_{\mathrm{t}}$
$:=$ 8 means that we will get 8 significance digits of the value
of $c_{2}$.

In the algorithm above, we have assumed implicitly that the average
error is a smooth and continuous function of $a$, $b$, $c_{2}$
for the values of $x_{k}=20k+100$, ($k$ = 1, 2, $\cdots$, 395).
For every $c_{2}$, we can get the value of $a$ and $b$, then obtain
the the average error $s_{1}$, so $s_{1}$ could be believed as a
convex and smooth function of $c_{2}$ (hence it will have only one
minimum point) in the interval we are considering. This could be verified
by plotting the figure of the curve $s_{1}$ = $s_{1}(c_{2})$ in
the given interval (although this work is not easy in practice).

If $n$ ranges from 120 to 8000 (step 20), we can get a fitting function
of $C_{1}(n)$, 
\begin{equation}
y=\dfrac{-0.02651010067}{\sqrt{x+4.8444724}}-0.3456324524,\label{eq:C1-Fit-Result-new}
\end{equation}
 with a minimal average error $2.446731760\times10^{-7}$. 

If $n$ ranges from 120 to 6000 (step 20), the fitting function of
$C_{1}(n)$ is, 
\begin{equation}
y=\dfrac{-0.02649625326}{\sqrt{x+4.7152127}}-0.3456327903,
\end{equation}
 with a minimal average error $2.279396699\times10^{-7}$. 

In the next section,   \eqref{eq:C1-Fit-Result-new} will be used
to estimate $C_{1}(n)$, i.e., 
\begin{equation}
C_{1}(n)\doteq\dfrac{-0.02651010067}{\sqrt{n+4.8444724}}-0.3456324524.\label{eq:C1-Fit-Result-new-B}
\end{equation}

\section{Fit the Denominator \label{sec:Fit-the-Denominator-C2(n)}}

By  \eqref{eq:Ram-Hardy-Rev-Form1} and  \eqref{eq:C1-Fit-Result-new-B},
we have 
\begin{equation}
C_{2}(n)\doteq\dfrac{\exp\left(\sqrt{\frac{2}{3}}\pi\sqrt{n+C_{1}(n)}\right)}{4\sqrt{3}p(n)}-n.\label{eq:C2-fml}
\end{equation}

\begin{figure}[h]
\centering\includegraphics[scale=0.36]{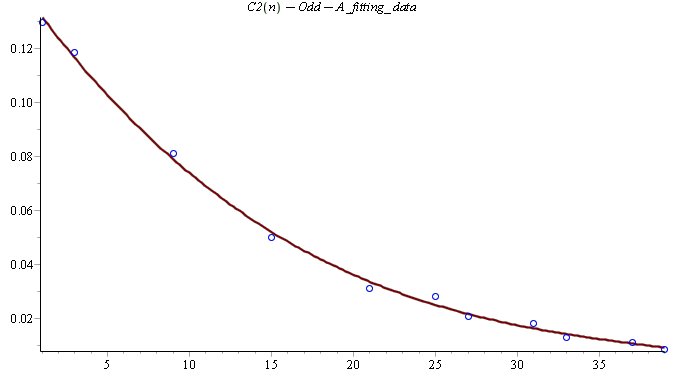}

\caption{Fit $\left(n,\,C_{2}(n)\right)$, the odd, Part A}
 \label{Fig:2_4_11-Fit-(n,C2)-Odd-A}

\hfill{}

\centering\includegraphics[scale=0.4]{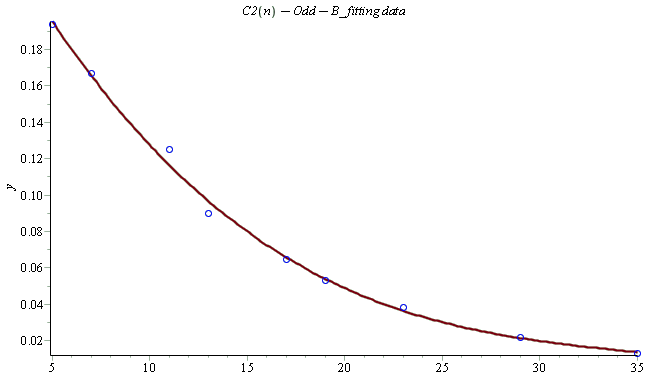}

\caption{Fit $\left(n,\,C_{2}(n)\right)$, the odd, Part B}
 \label{Fig:2_4_12-Fit-(n,C2)-Odd-B}
\end{figure}

If we point out the data $\left(n,\,C_{2}(n)\right)$ (1 $\leqslant$
$n$ $\leqslant$ 80) on the coordinate system as shown on  Figure
\ref{Fig:2_4_10-(n,C2)} on page \pageref{Fig:2_4_10-(n,C2)}, we
will immediately know than $C_{2}(n)$ can not be fit by a simple
function. From the  Figure \ref{Fig:2_4_10-(n,C2)} (or the value
of $C_{2}(n)$ calculated by a small program), it is clear that $C_{2}(n)$
is very small when $n$ $>$40, at least much less than $n$, so there
is no need to fit $C_{2}(n)$ when $n$ $>$ 40. 

When $n$ is odd, the points of $\left(n,\,C_{2}(n)\right)$ in Figure
\ref{Fig:2_4_10-(n,C2)} are above the horizontal-axis, it is not
difficult to separate them into two parts and fit them by two cubic
curves, as shown on Figure \ref{Fig:2_4_11-Fit-(n,C2)-Odd-A} and
Figure \ref{Fig:2_4_12-Fit-(n,C2)-Odd-B}. The two fitting functions
are 
\begin{align*}
y= & -1.548835311\times10^{-6}\times x^{3}+\\
 & \quad1.880663805\times10^{-4}\times x^{2}-\\
 & \quad0.008334098201\times x+0.1399798428,\\
y= & -5.416501948\times10^{-6}\times x^{3}+\\
 & \quad5.728510889\times10^{-4}\times x^{2}-\\
 & \quad0.02125835759\times x+0.2882706948.
\end{align*}

\begin{figure}[h]

\noindent \begin{centering}
\includegraphics[scale=0.26]{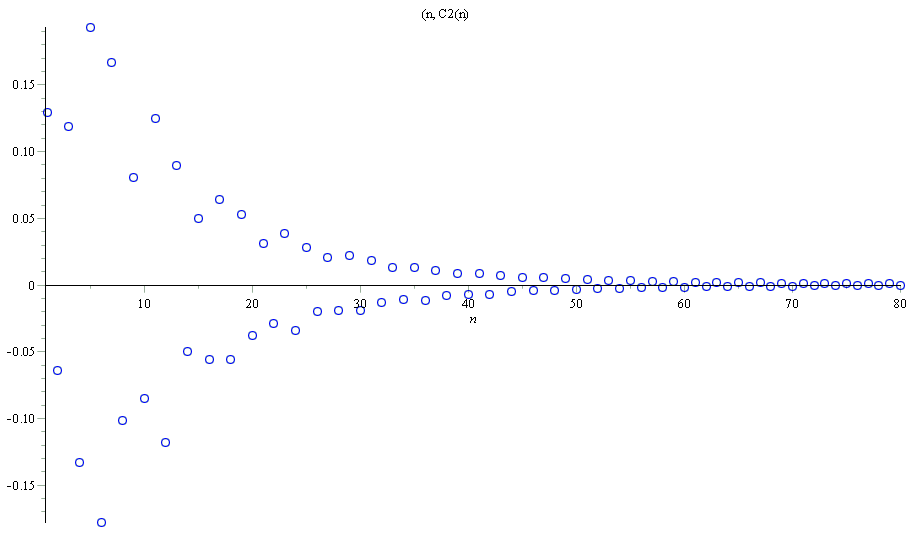}
\par\end{centering}

\caption{The graph of the data $\left(n,\,C_{2}(n)\right)$}
 \label{Fig:2_4_10-(n,C2)}\medskip{}

\noindent \begin{centering}
\includegraphics[scale=0.27]{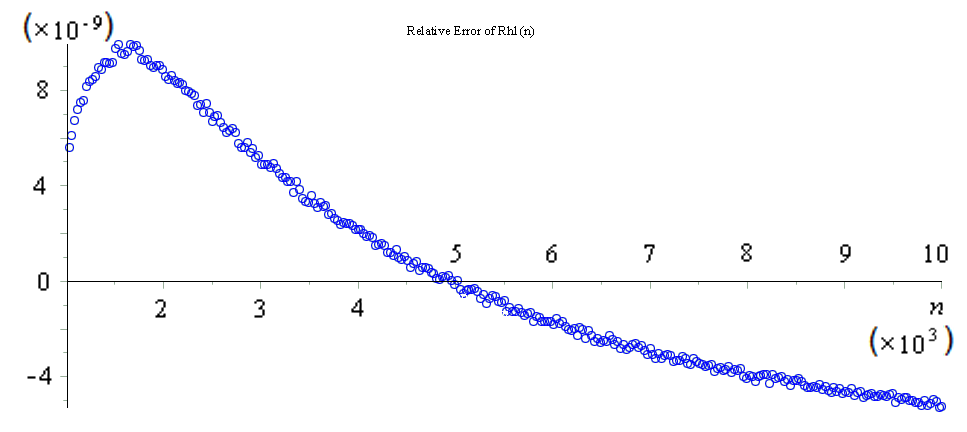}
\par\end{centering}

\caption{The Relative Error of $R_{\mathrm{h1}}(n)$ when 1000 $\leqslant$
$n$ $\leqslant$ 10000}
 \label{Fig:2_4_13-Rel-Err-Rh1}

\end{figure}

For the points of $\left(n,\,C_{2}(n)\right)$ under the horizontal-axis
(when $n$ is even) in Figure \ref{Fig:2_4_10-(n,C2)}, we have to
separate them into at least 4 parts so as to fit them smoothly, two
or three parts are not convenient.

As a result, we have to fit $C_{2}(n)$ by a hybrid function with
at least 6 pieces, or fit $p(n)$ by a piecewise-defined function
with 7 pieces, which is very complicated. This seems to contradict
with our purpose at the beginning of this paper.

From  Figure \ref{Fig:2_4_11-Fit-(n,C2)-Odd-A} on page \pageref{Fig:2_4_11-Fit-(n,C2)-Odd-A}
we found that the value of $C_{2}(n)$ are much less than $n$ when
$n$ $\geqslant$ 15, so the error will be very tiny if we omit $C_{2}(n)$.
 Hence we can calculate $p(n)$ directly by 
\begin{equation}
R_{\mathrm{h1}}(n)=\dfrac{1}{4\sqrt{3}n}\exp\left(\sqrt{\frac{2}{3}}\pi\sqrt{n+\dfrac{a_{1}}{\sqrt{n+c_{1}}}+b_{1}}\right),\label{eq:Ram-Hardy-Rev-Fml1A}
\end{equation}

where $a_{1}=-0.02651010067$, $b_{1}=-0.3456324524$ and $c_{1}=4.8444724$.

The error of  \eqref{eq:Ram-Hardy-Rev-Fml1A} to $p(n)$ (when $n$
$\leqslant$ 1000) is shown on Table \ref{Table:Rel-Err-p(n)-HR1(x)}
on page \pageref{Table:Rel-Err-p(n)-HR1(x)}. The accuracy is much
better than  \eqref{eq:Ram-Hardy-Fml-Est-1}. Although this fitting
function is obtained when $n$ $\geqslant$ 120, the relative error
is less than $6\times10^{-7}$ when $n$ $\geqslant$ 100, less than
$1\permil$ when $n$ $\geqslant$ 26, less than $1\%$ when $n$
$\geqslant$ 11. When 1000 $\leqslant$ $n$ $\leqslant$ 3000, the
relative error is less than $1\times10^{-8}$. When 3000 $\leqslant$
$n$ $\leqslant$ 10000, the relative error is less than $5.3\times10^{-9}$,
as shown on  Figure \ref{Fig:2_4_13-Rel-Err-Rh1} on page \pageref{Fig:2_4_13-Rel-Err-Rh1}.
But the relative error is not so satisfying when $n$ $\leqslant$
7, especially when $n=1$.

Consider that $p(n)$ is an integer, if we take the round approximation
of  \eqref{eq:Ram-Hardy-Rev-Fml1A}, \label{Sym:R'h1(n)} \nomenclature[Rh1(n)]{$R'_{\mathrm{h1}}(n)$}{The Hardy-Ramanujan's revised estimation formula 1. \pageref{Sym:R'h1(n)}}
\begin{equation}
R'_{\mathrm{h1}}(n)=\left\lfloor \dfrac{\exp\left(\sqrt{\frac{2}{3}}\pi\sqrt{n+\dfrac{a_{1}}{\sqrt{n+c_{1}}}+b_{1}}\right)}{4\sqrt{3}n}+\dfrac{1}{2}\right\rfloor ,\label{eq:Ram-Hardy-Rev-Fml-1}
\end{equation}
 (we may call it \emph{Hardy-Ramanujan's revised estimation formula
}1\index{Hardy-Ramanujan's revised estimation formula 1}), it will
solve perfectly the relative error problem when $n$ $<$ 11, as shown
on Table \ref{Table:Rel-Err-p(n)-HR1N(n)-round} on page \pageref{Table:Rel-Err-p(n)-HR1N(n)-round},
although the relative error will increase very little for some $n$,
which is negligible. (The average relative error is less than $2\times10^{-8}$
when $n$ $\geqslant$ 200.) Take an example, when $n=100$, $R_{\mathrm{h2}}(100)$
= 190569177, $p(100)$ = 190569292, the difference is 115; when $n=200$,
$R_{\mathrm{h2}}(200)$ = 3972999059745, $p(200)$ = 3972999029388,
the difference is 30357. Although the errors are much greater than
the error 0.004 of Hardy-Ramanujan formula with 6 terms ($n=100$)
or 8 terms ($n$ = 200) (refer \cite{Ramanujan1918AsymFmlComAnal}
or \cite{Rademacher1937ConvergSeries}), it contains only one term
of elementary functions, and is convenient for a junior middle school
student to calculate the value of $p(n)$ with high accuracy.

\begin{table}[h]
\noindent 

\noindent \begin{centering}
\includegraphics[bb=97bp 493bp 459bp 681bp,clip,scale=0.63]{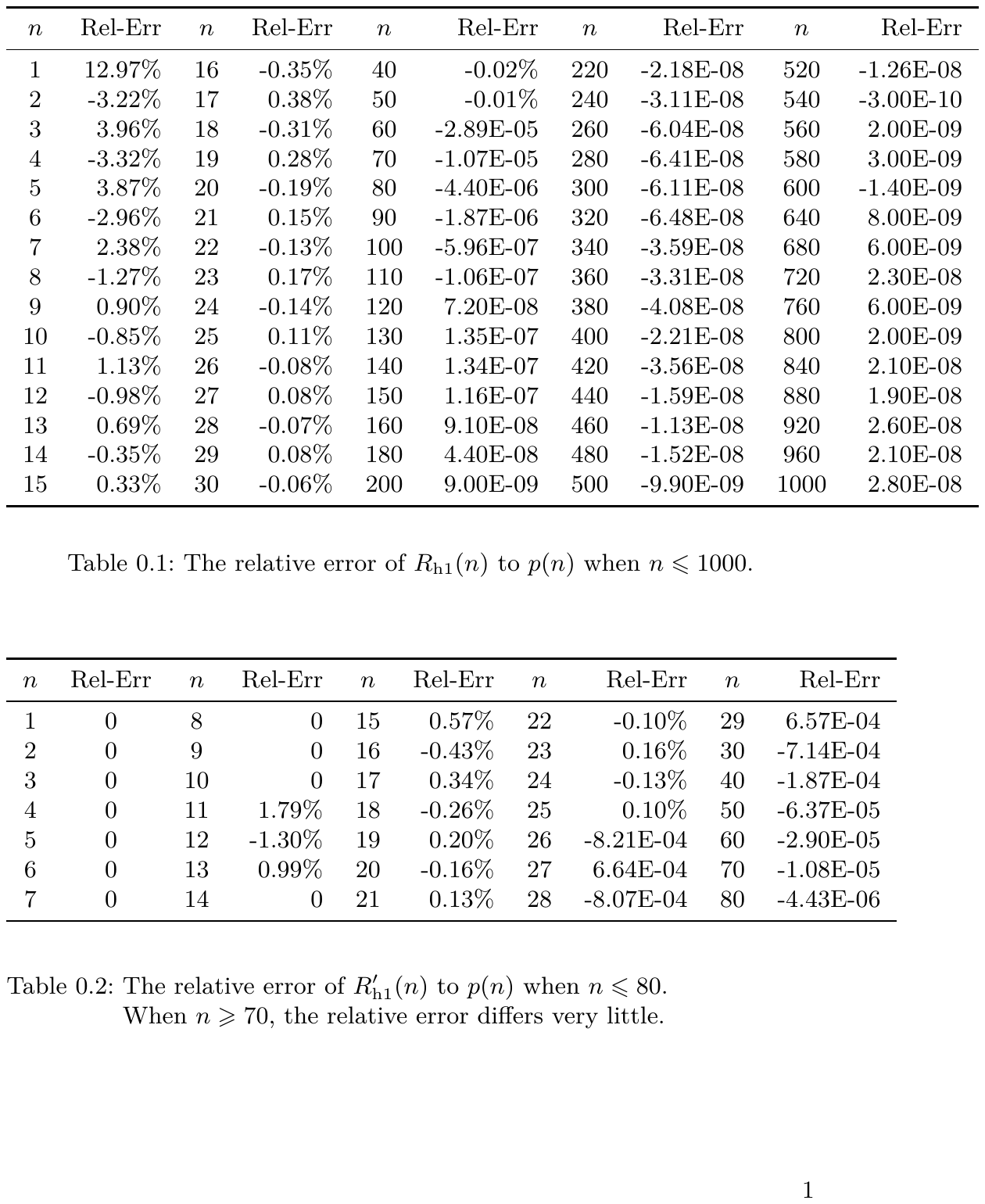}
\par\end{centering}

\caption{The relative error of $R_{\mathrm{h1}}(n)$ to $p(n)$ when $n\leqslant1000$.}
\label{Table:Rel-Err-p(n)-HR1(x)}

\vspace{0.8cm}

\noindent \begin{centering}
\includegraphics[bb=97bp 340bp 428bp 441bp,clip,scale=0.69]{tables-v3/Table-2_4_1_b-Rel-Err-RH1_n_-RH1N_n_round}
\par\end{centering}

\caption{The relative error of $R'_{\mathrm{h1}}(n)$ to $p(n)$ when $n\leqslant80$.\protect \\
When $n\geqslant70$, the relative error differs very little.}
\label{Table:Rel-Err-p(n)-HR1N(n)-round}

\end{table}

\section{Some Other Methods \label{sec:Estmt-p(n)-other-methd}}

In the previous sections, we assume that  $C_{1}\left(n\right)\doteq\frac{3\left(\ln\left(4n\sqrt{3}p(n)\right)\right)^{2}}{2\pi^{2}}-n$,
then fit the data $\left(n,\ \frac{3\left(\ln\left(4n\sqrt{3}p(n)\right)\right)^{2}}{2\pi^{2}}\right)$
($n=20k+100$, $k$ = 1, 2, $\cdots$, 395), and estimate $p(n)$
by $R_{\mathrm{h2}}(n)=\left\lfloor \frac{\exp\left(\sqrt{\frac{2}{3}}\pi\sqrt{n+C_{1}\left(n\right)}\right)}{4\sqrt{3}n}+\dfrac{1}{2}\right\rfloor $.

\subsection{Modify the Denominator only\label{sub:Modify the Denominator only}}

If we assume that $p(n)\doteq\frac{\exp\left(\pi\sqrt{\frac{2}{3}n}\right)}{4\sqrt{3}(n+C_{2})}$,
then 
\[
C_{2}(n)\doteq\dfrac{1}{4\sqrt{3}p(n)}\exp\left(\pi\sqrt{\frac{2}{3}n}\right)-n,
\]
 
\begin{figure}[h]

\noindent \begin{centering}
\includegraphics[scale=0.3]{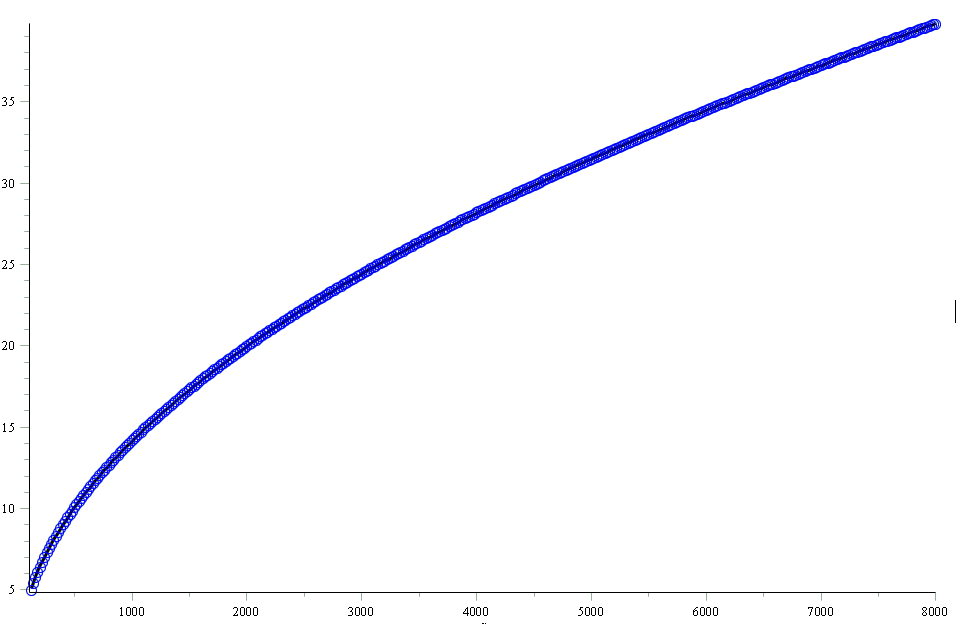}
\par\end{centering}

\noindent \begin{centering}
\caption{The graph of the data $\left(n,\ \frac{\exp\left(\pi\sqrt{\frac{2}{3}n}\right)}{4\sqrt{3}p(n)}-n\right)$
}
\label{Fig:2_4_14_C2(n)-Alone-no-C1-1}
\par\end{centering}

\end{figure}
we wonder whether we can fit the data $\left(n,\ \frac{\exp\left(\pi\sqrt{\frac{2}{3}n}\right)}{4\sqrt{3}p(n)}-n\right)$
 ($n=20k+100$, $k$ = 1, 2, $\cdots$, 395) by a function $C_{2}$
and estimate $p(n)$ by $\left\lfloor \frac{\exp\left(\pi\sqrt{\frac{2}{3}n}\right)}{4\sqrt{3}(n+C_{2})}+\dfrac{1}{2}\right\rfloor $?

The data $\left(n,\ \frac{\exp\left(\pi\sqrt{\frac{2}{3}n}\right)}{4\sqrt{3}p(n)}-n\right)$
 ($n=20k+100$, $k$ = 1, 2, $\cdots$, 395) are shown on  Figure
\ref{Fig:2_4_14_C2(n)-Alone-no-C1-1} on page \pageref{Fig:2_4_14_C2(n)-Alone-no-C1-1}
(together with the figure of a fitting function). It is not difficult
to know that a function in this form 
\[
y=a_{1}\times(x+c_{1})^{e_{1}}+b_{1}
\]
 will fit the points very well, and $e_{1}=0.5$ will be very satisfying.
By the same method to fit $C_{1}(n)$, we can obtain a fitting function
\[
y=0.4432884566\times\sqrt{x+0.274078}+0.1325096085
\]
 to fit $C_{2}(n)$ with an average error $3.65\times10^{-6}$.

Hence we can calculate $p(n)$ by 
\begin{equation}
R_{\mathrm{h2}}(n)=\dfrac{\exp\left(\sqrt{\frac{2}{3}}\pi\sqrt{n}\right)}{4\sqrt{3}\left(n+a_{2}\sqrt{n+c_{2}}+b_{2}\right)},\label{eq:Ram-Hardy-Rev-Fml2A}
\end{equation}

where $a_{2}=0.4432884566$, $b_{2}=0.1325096085$ and $c_{2}=0.274078$,
when $n$ is not so small.

\begin{table}[h]
\noindent 

\noindent \begin{centering}
\includegraphics[bb=97bp 489bp 459bp 676bp,clip,scale=0.63]{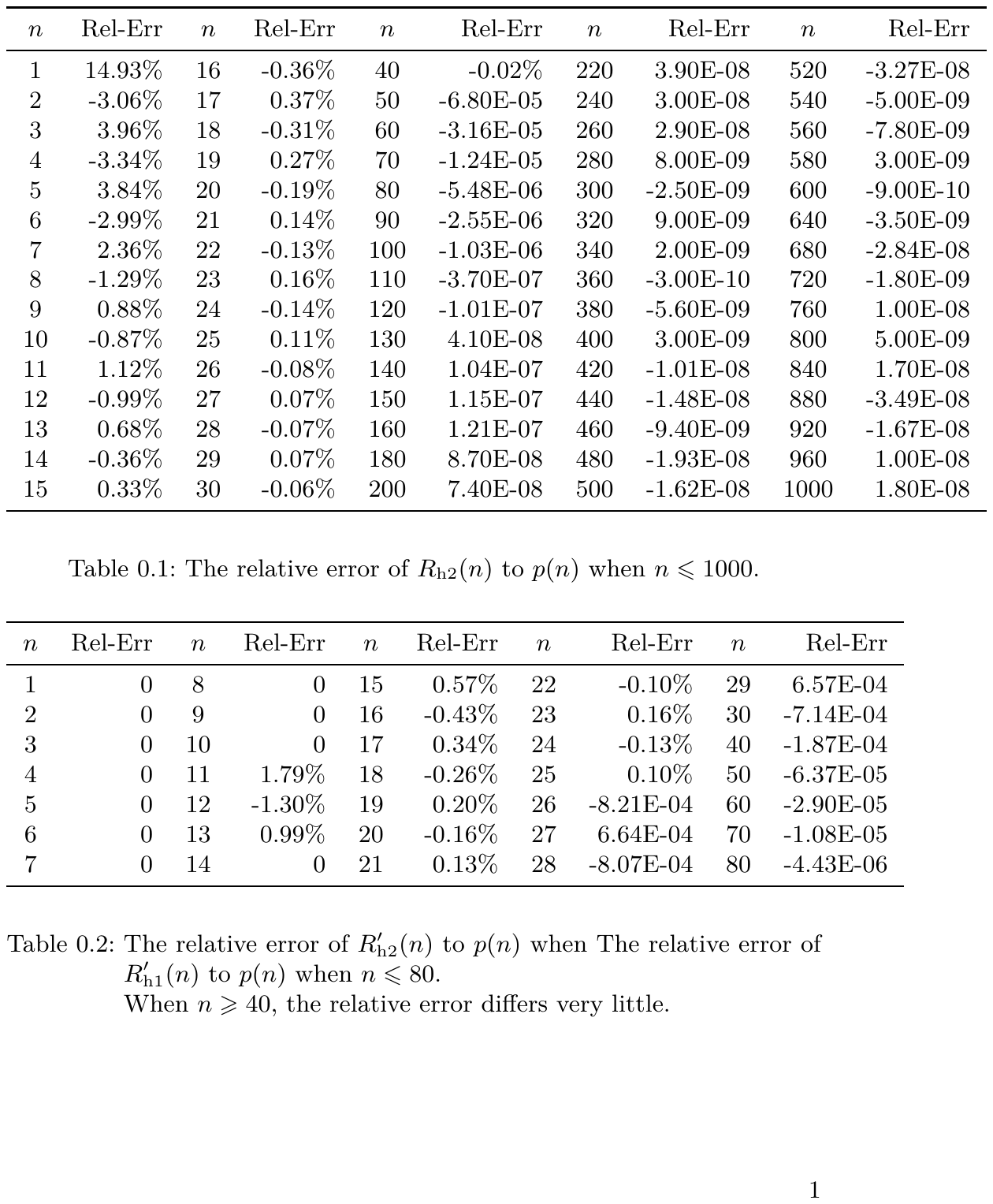}
\par\end{centering}

\caption{The relative error of $R_{\mathrm{h2}}(n)$ to $p(n)$ when $n\leqslant1000$.}
\label{Table:Rel-Err-p(n)-HR2(x)}

\vspace{0.8cm}

\noindent \begin{centering}
\includegraphics[bb=97bp 352bp 429bp 453bp,clip,scale=0.69]{tables-v3/Table-2_4_1_c-Rel-Err-RH2_n_-RH2N_n_round}
\par\end{centering}

\caption{The relative error of $R'_{\mathrm{h2}}(n)$ to $p(n)$ when $n\leqslant80$.\protect \\
When $n\geqslant40$, the relative error differs very little.}
\label{Table:Rel-Err-p(n)-HR2N(n)-round}

\end{table}

The error of  \eqref{eq:Ram-Hardy-Rev-Fml2A} to $p(n)$ is shown
on Table \ref{Table:Rel-Err-p(n)-HR2(x)} on page \pageref{Table:Rel-Err-p(n)-HR2(x)}
when $n$ $\leqslant$ 1000. The accuracy is much better than  \eqref{eq:Ram-Hardy-Fml-Est-1}.
Compared with Table \ref{Table:Rel-Err-p(n)-HR1(x)} (page \pageref{Table:Rel-Err-p(n)-HR1(x)}),
the accuracy are almost the same when $n$ $\leqslant$ 1000. When
1500 $\leqslant$ $n$ $\leqslant$ 10000, the relative error is obviously
less than that of  \eqref{eq:Ram-Hardy-Rev-Fml1A}, as shown on 
Figure \ref{Fig:2_4_15-Rel-Err-Rh2} on page \pageref{Fig:2_4_15-Rel-Err-Rh2}
(compared with  Figure \ref{Fig:2_4_14_C2(n)-Alone-no-C1-1} on page
\pageref{Fig:2_4_14_C2(n)-Alone-no-C1-1}). Which means that $R_{\mathrm{h2}}(n)$
is more accurate than $R_{\mathrm{h1}}(n)$. (If we change the range
of $n$ of the data points, the accuracy of the fitting function obtained
may not be so good.)

Consider that $p(n)$ is an integer, we can take the round approximation
of  \eqref{eq:Ram-Hardy-Rev-Fml2A}, \label{Sym:R'h2(n)} \nomenclature[Rh2(n)]{$R'_{\mathrm{h2}}(n)$}{The Hardy-Ramanujan's revised estimation formula 2. \pageref{Sym:R'h2(n)}}
\begin{equation}
R'_{\mathrm{h2}}(n)=\left\lfloor \dfrac{\exp\left(\sqrt{\frac{2}{3}}\pi\sqrt{n}\right)}{4\sqrt{3}\left(n+a_{2}\sqrt{n+c_{2}}+b_{2}\right)}+\dfrac{1}{2}\right\rfloor ,\label{eq:Ram-Hardy-Rev-Fml-2}
\end{equation}
 for small values of $n$. We may call it \emph{Hardy-Ramanujan's
revised estimation formula }2. \index{Hardy-Ramanujan's revised estimation formula 2}
The error of  \eqref{eq:Ram-Hardy-Rev-Fml-2} to $p(n)$ is shown
on Table \ref{Table:Rel-Err-p(n)-HR2N(n)-round} (on page \pageref{Table:Rel-Err-p(n)-HR2N(n)-round})
when $n$ $\leqslant$ 1000.

\subsection{Fit $\nicefrac{R_{\mathrm{h}}(n)}{p(n)}$}

At the beginning of section \ref{sec:main-idea} , some other methods
to estimate $p(n)$ are mentioned, such as estimating the value of
$\frac{R_{\mathrm{h}}(n)}{p(n)}$ by a function $f_{1}(n)$, then
estimate $p(n)$ by $\frac{R_{\mathrm{h}}(n)}{f_{1}(n)}$. 

The data $\left(n,\ \frac{R_{\mathrm{h}}(n)}{f_{1}(n)}\right)$ ($n=20k+100$,
$k$ = 1, 2, $\cdots$, 395) are shown on  Figure \ref{Fig:2_4_16_Rh-over-p(n)}
on page \pageref{Fig:2_4_16_Rh-over-p(n)} (together with the figure
of a fitting function). It is not difficult to find out that a function
\[
y=1+\dfrac{1}{\sqrt{a_{3}x+b_{3}}},
\]
where $a_{3}=5.062307637$ and $b_{3}=-75.65700620$, will fit the
data very well, as shown on the figure, with an average error $1.41\times10^{-4}$.
(because the data $\left(n,\ \left(\frac{R_{\mathrm{h}}(n)}{f_{1}(n)}-1\right)^{-2}\right)$
lies exactly on a straight line $y=a_{3}x+b_{3}$, as shown on Figure
\ref{Fig:2_4_17-RH(n)-over-p(n)-mid-data-fit} on page \pageref{Fig:2_4_17-RH(n)-over-p(n)-mid-data-fit})

\begin{figure}[h]
\centering \includegraphics[scale=0.3]{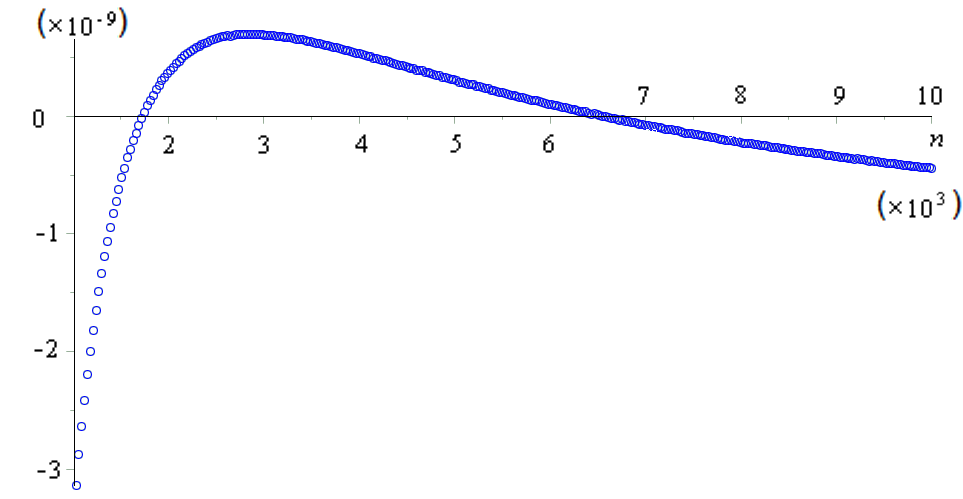}

\caption{The Relative Error of $R_{\mathrm{h2}}(n)$ when 1000 $\leqslant$
$n$ $\leqslant$ 10000}
 \label{Fig:2_4_15-Rel-Err-Rh2}

$\ $\medskip{}

\centering\includegraphics[scale=0.3]{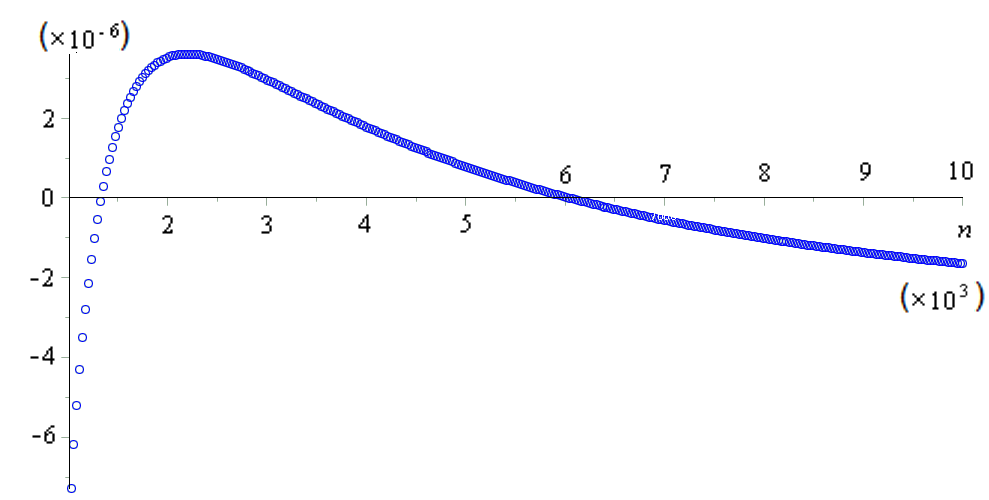}

\caption{The Relative Error of $R_{d\mathrm{3}}(n)$ when 1000 $\leqslant$
$n$ $\leqslant$ 10000}
 \label{Fig:2_4_18-Rel-Err-R3}
\end{figure}

So we have another fitting function for $p(n)$, 
\[
R_{d3}(n)=\dfrac{R_{\mathrm{h}}(n)}{1+\dfrac{1}{\sqrt{a_{3}n+b_{3}}}}.
\]

However, this formula does not fit $p(n)$ very well when $n$ is
small. When $n\leqslant14$, the value of $R_{d3}(n)$ is an imaginary
number. Unfortunately, when $n$ > 1000, the error of $R_{d3}(n)$
to $p(n)$ is about 1000 times of the error of $R_{\mathrm{h2}}(n)$,
as shown on  Figure \ref{Fig:2_4_18-Rel-Err-R3} on page \pageref{Fig:2_4_18-Rel-Err-R3}.

Actually, $R_{\mathrm{h2}}(n)$ is in the form $\frac{R_{\mathrm{h}}(n)}{f_{1}(n)}$,
since $\frac{\exp\left(\sqrt{\frac{2}{3}}\pi\sqrt{n}\right)}{4\sqrt{3}\left(n+a_{2}\sqrt{n+c_{2}}+b_{2}\right)}$
= $\frac{\exp\left(\sqrt{\frac{2}{3}}\pi\sqrt{n}\right)}{4\sqrt{3}n}$
$\frac{n}{n+a_{2}\sqrt{n+c_{2}}+b_{2}}$ = $\frac{R_{\mathrm{h}}(n)}{1+\frac{a_{2}}{n}\sqrt{n+c_{2}}+\frac{b_{2}}{n}}$.
As $1+\frac{a_{2}\sqrt{n+c_{2}}}{n}+\frac{b_{2}}{n}$ fits $\frac{R_{\mathrm{h}}(n)}{p(n)}$
with very little error, $1+\frac{1}{\sqrt{a_{3}n+b3}}$ will not reach
that accuracy.

\begin{figure}[h]
\centering\includegraphics[scale=0.26]{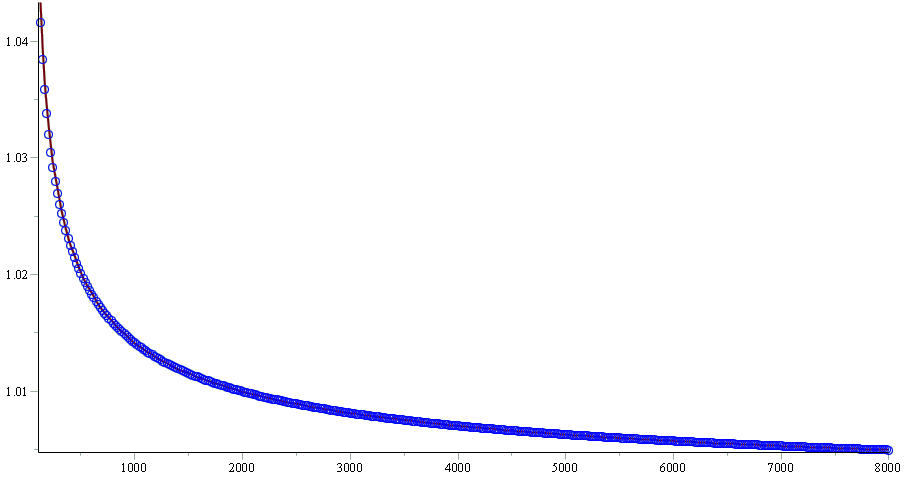}

\caption{The graph of the data $\left(n,\ \frac{R_{\mathrm{h}}(n)}{p(n)}\right)$
 and the fitting function}
\label{Fig:2_4_16_Rh-over-p(n)}

$\ $\medskip{}

\centering\includegraphics[scale=0.25]{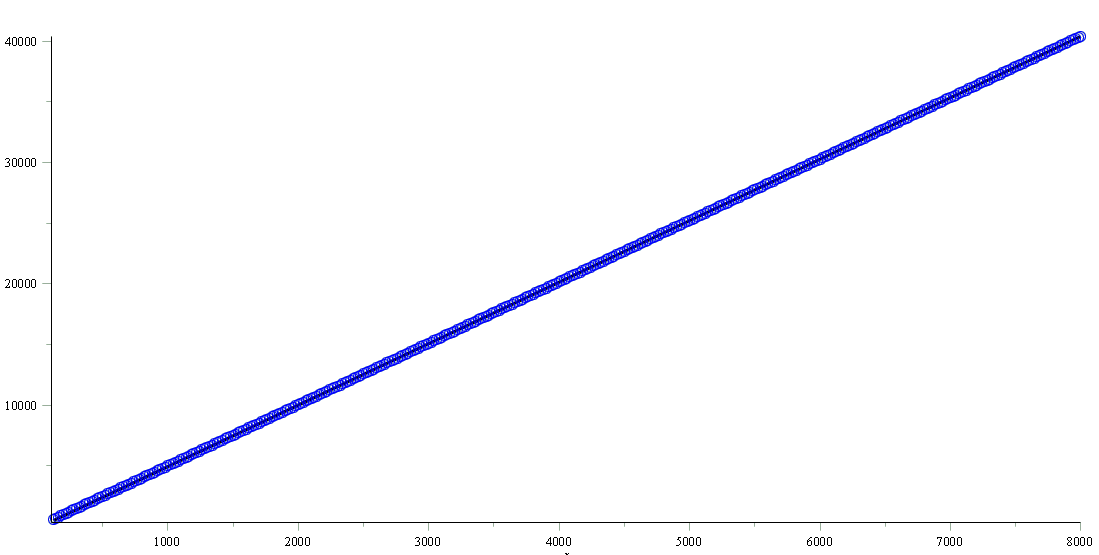}

\caption{The data $\left(n,\ \left(\frac{R_{\mathrm{h}}(n)}{f_{1}(n)}-1\right)^{-2}\right)$
and the fitting function}
 \label{Fig:2_4_17-RH(n)-over-p(n)-mid-data-fit}
\end{figure}

\section{Approximate $p(n)$ by Fitting $R_{\mathrm{h}}(n)-p(n)$ \label{sec:Estimate-p(n)-by-fitting-Rh(n)-p(n)}}

\subsection{Result 1}

It is not difficult to verify that 
\begin{equation}
R_{\mathrm{h}}(n)-R_{\mathrm{h}}(n-1)\sim\dfrac{\pi}{12\sqrt{2n^{3}}}\exp\left(\sqrt{\frac{2}{3}}\pi\sqrt{n}\right).\label{eq:rh(n)-rh(n-1)_est}
\end{equation}

(refer  Sec. 3.1 of \cite{liwenwei2016-Num-Cnj-Cls-Der-arXiv}).
As $R_{\mathrm{h}}(n)$ is obviously greater than $p(n)$, we wander
whether we can fit $R_{\mathrm{h}}(n)-p(n)$ by an expression similar
like the right part of \ref{eq:rh(n)-rh(n-1)_est}, such as $\frac{\pi\exp\left(\sqrt{\frac{2}{3}}\pi\sqrt{n}\right)}{12\sqrt{2C_{3}(n)}}$,
where $C_{3}(n)$ is a cubic function, or equivalently, fit $\left(\frac{\pi\exp\left(\sqrt{\frac{2}{3}}\pi\sqrt{n}\right)}{12\sqrt{2}\left(R_{\mathrm{h}}(n)-p(n)\right)}\right)^{2}$
by a cubic function $C_{3}(n)$, from the data with the data $(n,\,p(n))$
($n=20k+60$, $k$ = 1, 2, $\cdots$, 397).  The result is 
\[
C_{3}(n)=a_{1}n^{3}+b_{1}n^{2}+c_{1}n+d_{1},
\]
 where 
\begin{align*}
a_{1} & =8.383485427,\\
b_{1} & =130.0792015,\\
c_{1} & =-1.197477259\times10^{5},\\
d_{1} & =4.188653689\times10^{7}.
\end{align*}

\begin{figure}
\centering\includegraphics[scale=0.26]{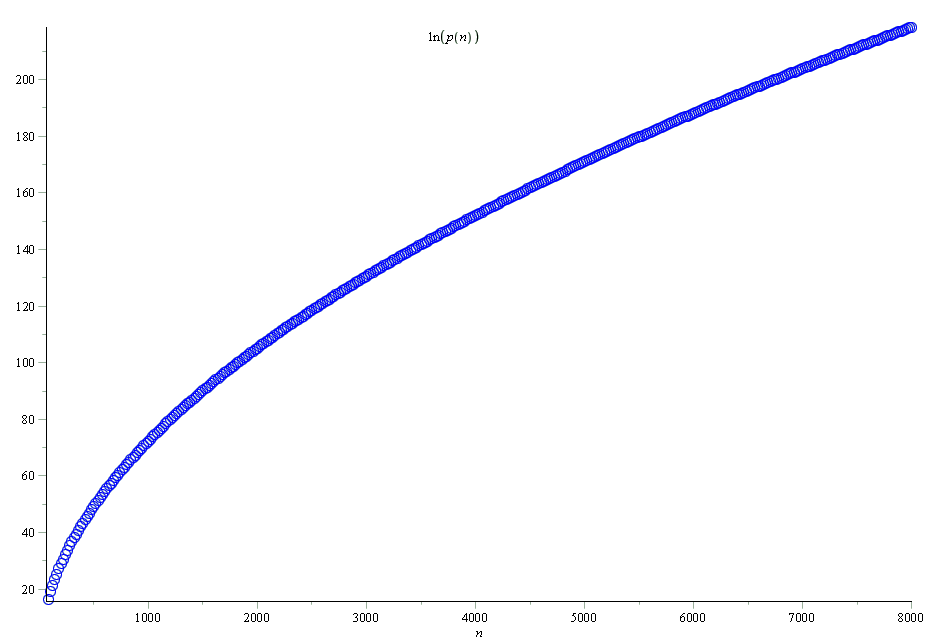}

\caption{The graph of the data $\bigl(n,\,\ln\left(p\left(n\right)\right)\bigr)$}
\label{Fig:2_4_19_ln(p(n))}

$\ $\medskip{}

\centering\includegraphics[scale=0.26]{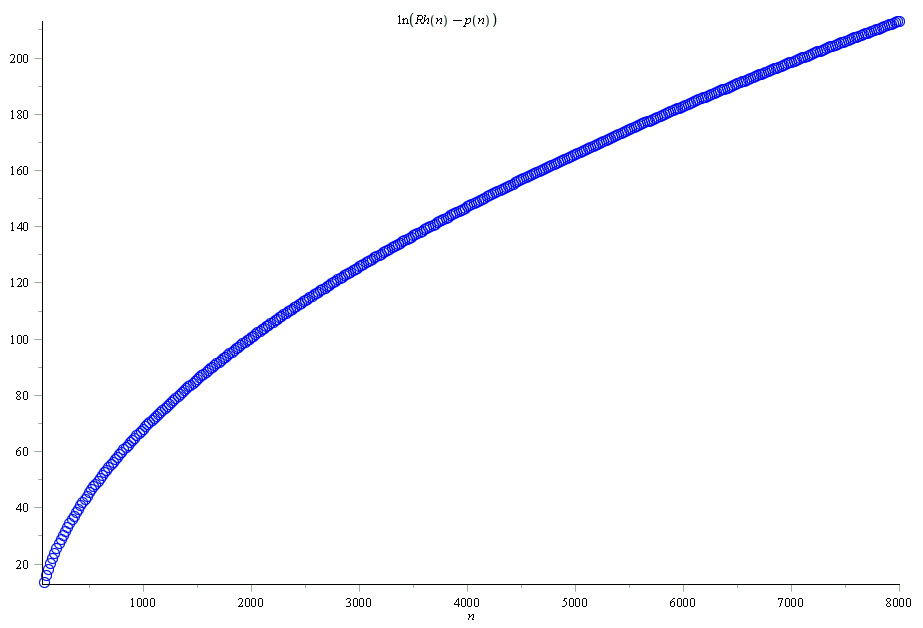}

\caption{The graph of the data $\bigl(n,\,\ln\left(R_{\mathrm{h}}(n)-p\left(n\right)\right)\bigr)$}
\label{Fig:2_4_20_ln(Rh(n)-p(n))}
\end{figure}

Here $c_{1}$ and $d_{1}$ are very huge, which suggests that this
result may not be so satisfying. As a sequence, if we estimate $p(n)$
by 
\[
F_{3}(n)=R_{\mathrm{h}}(n)-\dfrac{\pi\exp\left(\sqrt{\frac{2}{3}}\pi\sqrt{n}\right)}{12\sqrt{2C_{3}(n)}},
\]
 the relative error differs very little with the relative error of
$R_{\mathrm{h}}(n)$ to $p(n)$ when $n<50$, but the relative error
is not satisfying when $n<280$, as shown in Table \ref{Table:Rel-Err-p(n)-F3(n)}
on page \pageref{Table:Rel-Err-p(n)-F3(n)}.

If we fit $\left(\frac{\pi\exp\left(\sqrt{\frac{2}{3}}\pi\sqrt{n}\right)}{12\sqrt{2}\left(R_{\mathrm{h}}(n)-p(n)\right)}\right)^{2}$
by a function like 
\[
C_{3}(n)=a_{2}n^{3}+b_{2}n^{2.5}+c_{2}n^{2}+d_{2}n^{1.5}+e_{2}n+f_{2}n^{0.5}+g_{2},
\]
 the result are even worse, since imaginary number appeared (as concerned
to the data mentioned in this section. If we fit less data, the imaginary
problem might be avoid).

So we have to consider a different method.

\subsection{Result 2\label{sub:Result-2}}

In the previous sub-subsection, we obtained the asymptotic order of
$p(n)-p(n-1)$, and revised it to fit $R_{\mathrm{h}}(n)-p(n)$. Since
$R_{\mathrm{h}}(n)$ is always a little greater than $p(n)$, we may
guess that there is a $t_{0}$ such that $R_{\mathrm{h}}(n-t_{0})$
is closer to $p(n)$ than $R_{\mathrm{h}}(n)$. Then we can find the
asymptotic order of $R_{\mathrm{h}}(n)-R_{\mathrm{h}}(n-t_{0})$ and
use the new asymptotic order to fit $R_{\mathrm{h}}(n)-p(n)$.

\begin{table}
\noindent 

\noindent \begin{centering}
\includegraphics[bb=97bp 409bp 459bp 597bp,clip,scale=0.63]{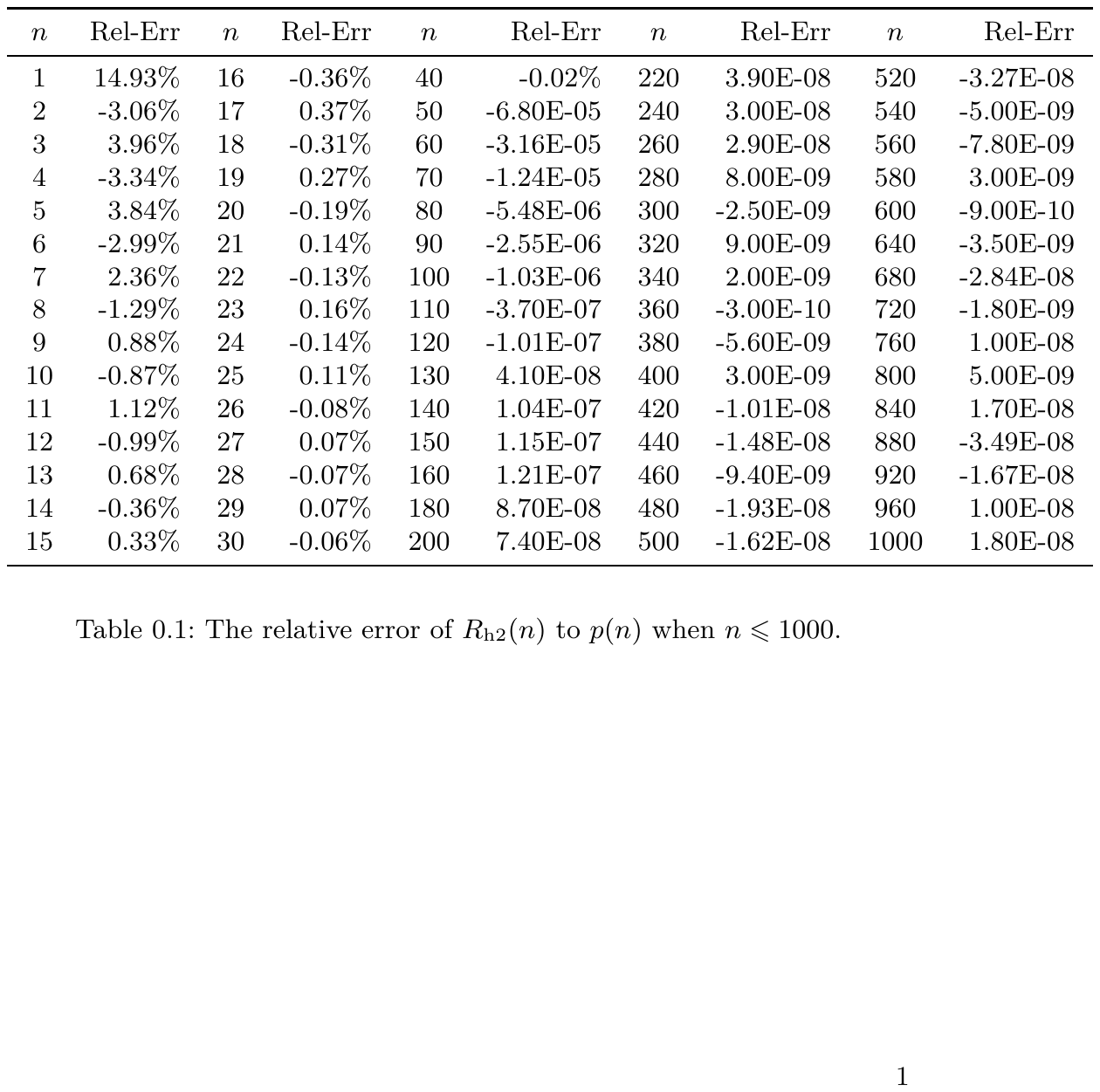}
\par\end{centering}

\caption{The relative error of $F_{3}(n)$ to $p(n)$ when $n\leqslant1000$.}
\label{Table:Rel-Err-p(n)-F3(n)}

\end{table}

By the same idea described in the algorithm mentioned on page \pageref{Alg:gen-c2-by-loop},
we can obtain the value $t_{0}\doteq0.3594143172$. 

When $n\gg1$ and $n\gg t$, 

$\ $ \hphantom{a } $r(n)$ $=$ $R_{\mathrm{h}}(n)-R_{\mathrm{h}}(n-t)$

$=$ $\frac{\exp\left(\sqrt{\frac{2}{3}}\pi\sqrt{n}\right)}{4\sqrt{3}n}-$
$\frac{\exp\left(\sqrt{\frac{2}{3}}\pi\sqrt{n-t}\right)}{4\sqrt{3}(n-t)}$

=  $\frac{\exp\left(\sqrt{\frac{2}{3}}\pi\sqrt{n}\right)}{4\sqrt{3}n}-\left(\frac{\exp\left(\frac{t\pi\sqrt{2/3}}{\sqrt{n}+\sqrt{n-t}}\right)}{n}-\frac{1}{(n-t)}\right)$

$\sim$ $\frac{\exp\left(\sqrt{\frac{2}{3}}\pi\sqrt{n-t}\right)}{4\sqrt{3}n}\left(\frac{\exp\left(\frac{t\pi\sqrt{2/3}}{2\sqrt{n-t/2}}\right)}{n}-\frac{1}{(n-t)}\right)$

$\sim$ $\frac{\exp\left(\sqrt{\frac{2}{3}}\pi\sqrt{n-t}\right)}{4\sqrt{3}n}\left(\frac{1+\frac{t\pi}{\sqrt{6(n-t/2)}}}{n}-\frac{1}{(n-t)}\right)$

$\sim$ $\frac{\exp\left(\sqrt{\frac{2}{3}}\pi\sqrt{n-t}\right)}{4\sqrt{3}n}\left(\frac{t\pi\sqrt{n+t/2}}{\sqrt{6}n(n-t)}\right)$

= $\frac{t\pi\exp\left(\sqrt{\frac{2}{3}}\pi\sqrt{n-t}\right)}{12\sqrt{2}(n-t)\sqrt{(n-t/2)}}$ 

$\sim$ $\frac{t\pi}{12\sqrt{2}\sqrt{n^{3}}}\exp\left(\sqrt{\frac{2}{3}}\pi\sqrt{n}\right).$

As 
\begin{equation}
r(n)\sim\dfrac{t\pi}{12\sqrt{2}\sqrt{n^{3}}}\exp\left(\sqrt{\frac{2}{3}}\pi\sqrt{n}\right),\label{eq:Deduction}
\end{equation}

so we may consider fitting $R_{\mathrm{h}}(n)-p(n)$ by $\frac{\sqrt{2}t_{0}\pi\exp\left(\sqrt{\frac{2}{3}}\pi\sqrt{n-t_{0}}\right)}{24C_{4}(n)}$,
where 
\begin{equation}
C_{4}(n)=a_{2}(n-t_{0})^{1.5}+b_{2}(n-t_{0})+c_{2}(n-t_{0})^{0.5}+d_{2}.\label{eq:C4(n)}
\end{equation}

\begin{table}[h]
\noindent 

\noindent \begin{centering}
\includegraphics[bb=97bp 479bp 459bp 666bp,clip,scale=0.63]{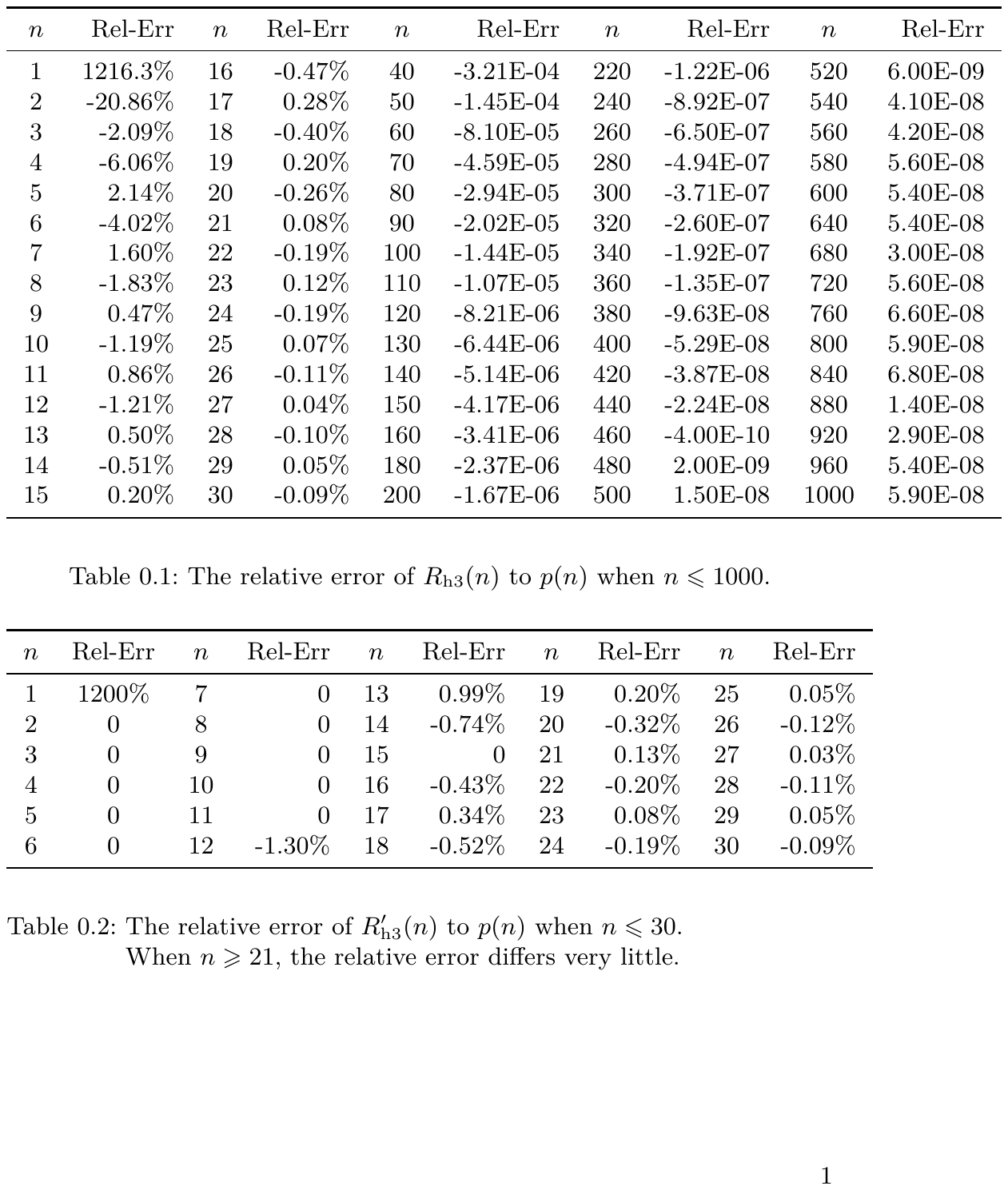}
\par\end{centering}

\caption{The relative error of $R_{\mathrm{h3}}(n)$ to $p(n)$ when $n\leqslant1000$.}
\label{Table:Rel-Err-p(n)-Rh3(n)}

\vspace{0.8cm}

\noindent \begin{centering}
\includegraphics[bb=97bp 352bp 411bp 442bp,clip,scale=0.73]{tables-v3/Table-2_4_1_e-Rel-Err-Rh3_n_-Rh3_n_round}
\par\end{centering}

\caption{The relative error of $R'_{\mathrm{h3}}(n)$ to $p(n)$ when $n\leqslant30$.\protect \\
When $n\geqslant21$, the relative error differs very little.}
\label{Table:Rel-Err-p(n)-Rh3(n)-round}

\end{table}

When $t_{0}\doteq0.3594143172$, \footnote{$\ $ In \cite{Sandor1995UnstrctPttInt} (or \cite{Sandor2006UnstrctPttInt})
or some other papers, there is a theoretic value $\dfrac{1}{24}$. } it is not difficult to find out that 
\begin{align*}
a_{2} & =1.039888529,\\
b_{2} & =-0.3305606395,\\
c_{2} & =0.6134039843,\\
d_{2} & =-0.8582793693,
\end{align*}
 from the data $(n,\,p(n))$ ($n=20k+60$, $k$ = 1, 2, $\cdots$,
397). Here none of the coefficients is very huge, which seems better
than the previous result mentioned in this section.  As a matter
of fact, if we estimate $p(n)$ by 
\begin{equation}
R_{\mathrm{h3}}(n)=R_{\mathrm{h}}(n)-\dfrac{\sqrt{2}t_{0}\pi\exp\left(\sqrt{\frac{2}{3}}\pi\sqrt{n-t_{0}}\right)}{24C_{4}(n)},\label{eq:Ram-Hardy-Rev-Fml-3A}
\end{equation}
 the relative error is very small even when $n<10$ (except the cases
when $n$ = 1 or 2) as shown on Table \ref{Table:Rel-Err-p(n)-Rh3(n)}
on page \pageref{Table:Rel-Err-p(n)-Rh3(n)}. This is the first time
to have an estimation formula of $p(n)$ which can reach a good accuracy
without taking round approximation even when $n<10$. 

Further more, if we take the round value of $R_{\mathrm{h3}}(n)$,
\label{Sym:R'h3(n)} \nomenclature[Rh3(n)]{$R'_{\mathrm{h3}}(n)$}{The Hardy-Ramanujan's revised estimation formula 3. \pageref{Sym:R'h3(n)}}
\begin{equation}
R'_{\mathrm{h3}}(n)=\left\lfloor R_{\mathrm{h}}(n)-\dfrac{\sqrt{2}t_{0}\pi\exp\left(\sqrt{\frac{2}{3}}\pi\sqrt{n-t_{0}}\right)}{24C_{4}(n)}+\dfrac{1}{2}\right\rfloor ,\label{eq:Ram-Hardy-Rev-Fml-3}
\end{equation}
the relative error to error is even less, especially when $n=15$
or $1<n<12$ (it reaches 0), as shown on Table \ref{Table:Rel-Err-p(n)-Rh3(n)-round}
on page \pageref{Table:Rel-Err-p(n)-Rh3(n)-round}. The relative error
is less than $3\times10^{-9}$ when $2500<n<10000$, as shown on 
Figure \ref{Fig:2_4_19_ln(p(n))} on page \pageref{Fig:2_4_19_ln(p(n))}.
This formula will be called \emph{Hardy-Ramanujan's revised estimation
formula }3. \index{Hardy-Ramanujan's revised estimation formula 3} 

\begin{figure}
\centering\includegraphics[scale=0.25]{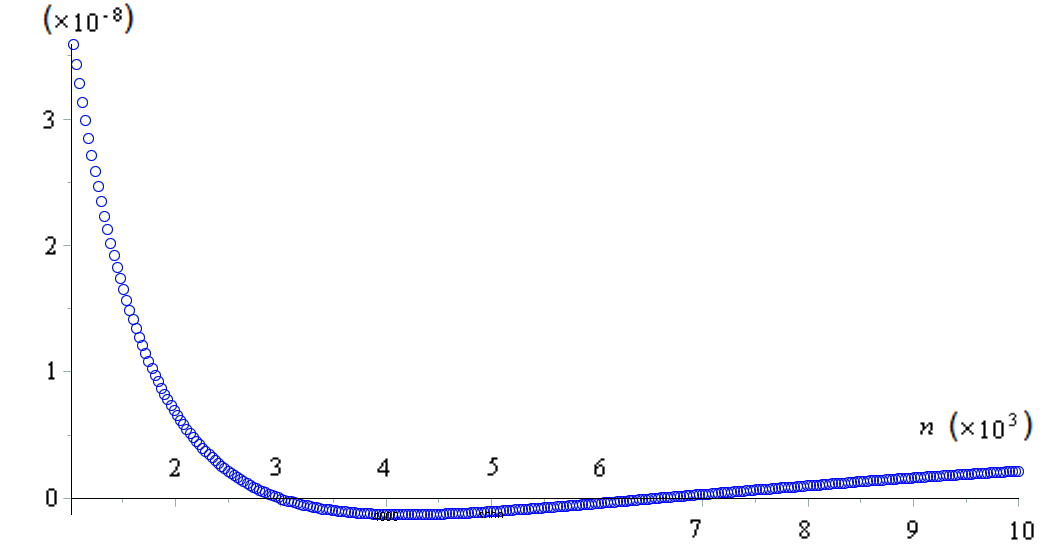}

\caption{The Relative Error of $R_{\mathrm{\mathrm{h}3}}(n)$ when 1000 $\leqslant$
$n$ $\leqslant$ 10000}
 \label{Fig:2_4_21-Rel-Err-Rh3}

$\ $\medskip{}

\centering\includegraphics[scale=0.25]{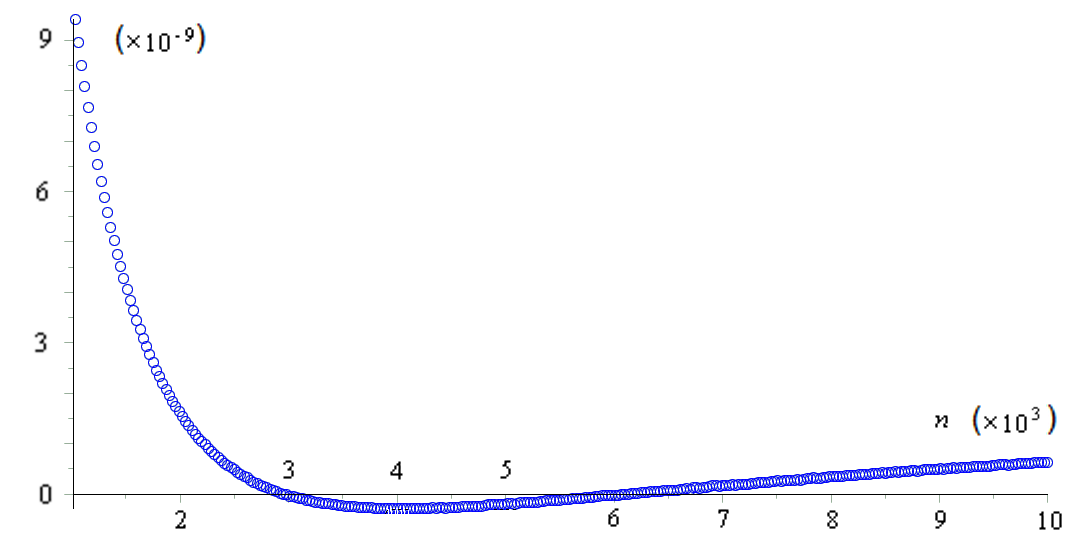}

\caption{The Relative Error of $R_{\mathrm{\mathrm{h}4}}(n)$ when 1000 $\leqslant$
$n$ $\leqslant$ 10000}
 \label{Fig:2_4_22-Rel-Err-Rh4}
\end{figure}

\subsection{Result 3}

Now that we can fit $R_{\mathrm{h}}(n)-p(n)$ by $\frac{\sqrt{2}t_{0}\pi\exp\left(\sqrt{\frac{2}{3}}\pi\sqrt{n-t_{0}}\right)}{24C_{4}(n)}$,
where\\
 $C_{4}(n)=a_{2}(n-t_{0})^{1.5}+b_{2}(n-t_{0})+c_{2}(n-t_{0})^{0.5}+d_{2}$,
maybe we can also fit $R_{\mathrm{h}}(n)-p(n)$ by $\frac{\pi\exp\left(\sqrt{\frac{2}{3}}\pi\sqrt{n}\right)}{12\sqrt{2}C_{5}(n)}$
directly, where 
\begin{equation}
C_{5}(n)=a_{3}n^{1.5}+b_{3}n+c_{3}n^{0.5}+d_{3},\label{eq:C5(n)}
\end{equation}
 or equivalently, to fit $\frac{\pi\exp\left(\sqrt{\frac{2}{3}}\pi\sqrt{n}\right)}{12\sqrt{2}\left(R_{\mathrm{h}}(n)-p(n)\right)}$
by a function $C_{5}(n)$ in the form mentioned above.

We can easily obtain the value of the unknown coefficients in the
equation above by the least square method. 
\begin{align*}
a_{3} & =2.893270736,\\
b_{3} & =0.4164546941,\\
c_{3} & =-0.08501098214,\\
d_{3} & =-0.4621004962.
\end{align*}
 Again, none of the coefficients is very huge. As a result, the relative
error of 
\begin{equation}
R_{\mathrm{h4}}(n)=R_{\mathrm{h}}(n)-\dfrac{\pi\exp\left(\sqrt{\frac{2}{3}}\pi\sqrt{n}\right)}{12\sqrt{2}C_{5}(n)},\label{eq:Ram-Hardy-Rev-Fml-4A}
\end{equation}
 to $p(n)$ is also very small when $n<10$ (even in the cases when
$n$ = 1 or 2) as shown on Table \ref{Table:Rel-Err-p(n)-Rh4(n)}
on page \pageref{Table:Rel-Err-p(n)-Rh4(n)}. This is the first time
to obtain an estimation formula of $p(n)$ which can reach a very
good accuracy even when $n<10$.

Further more, if we get the round value of $R_{\mathrm{h4}}(n)$,
\label{Sym:R'h4(n)} \nomenclature[Rh4(n)]{$R'_{\mathrm{h4}}(n)$}{The Hardy-Ramanujan's revised estimation formula 4. \pageref{Sym:R'h4(n)}}
\begin{equation}
R'_{\mathrm{h4}}(n)=\left\lfloor R_{\mathrm{h4}}(n)-\dfrac{\pi\exp\left(\sqrt{\frac{2}{3}}\pi\sqrt{n}\right)}{12\sqrt{2}C_{5}(n)}+\dfrac{1}{2}\right\rfloor ,\label{eq:Ram-Hardy-Rev-Fml-4}
\end{equation}
the relative error to error is even less, especially when $n=15$
or $1<n<12$ it reaches 0, as shown on Table \ref{Table:Rel-Err-p(n)-Rh4(n)-round}
on page \pageref{Table:Rel-Err-p(n)-Rh4(n)-round}. The relative error
is less than $1\times10^{-9}$ when $2500<n<10000$, as shown on 
Figure \ref{Fig:2_4_20_ln(Rh(n)-p(n))} on page \pageref{Fig:2_4_20_ln(Rh(n)-p(n))}.
That is much better than $R_{\mathrm{h3}}(n)$ and $R'_{\mathrm{h3}}(n)$,
besides, it is more simple. This formula will be called \emph{Hardy-Ramanujan's
revised estimation formula }4. \index{Hardy-Ramanujan's revised estimation formula 4} 

\begin{table}[h]
\noindent 

\noindent \begin{centering}
\includegraphics[bb=97bp 478bp 459bp 665bp,clip,scale=0.63]{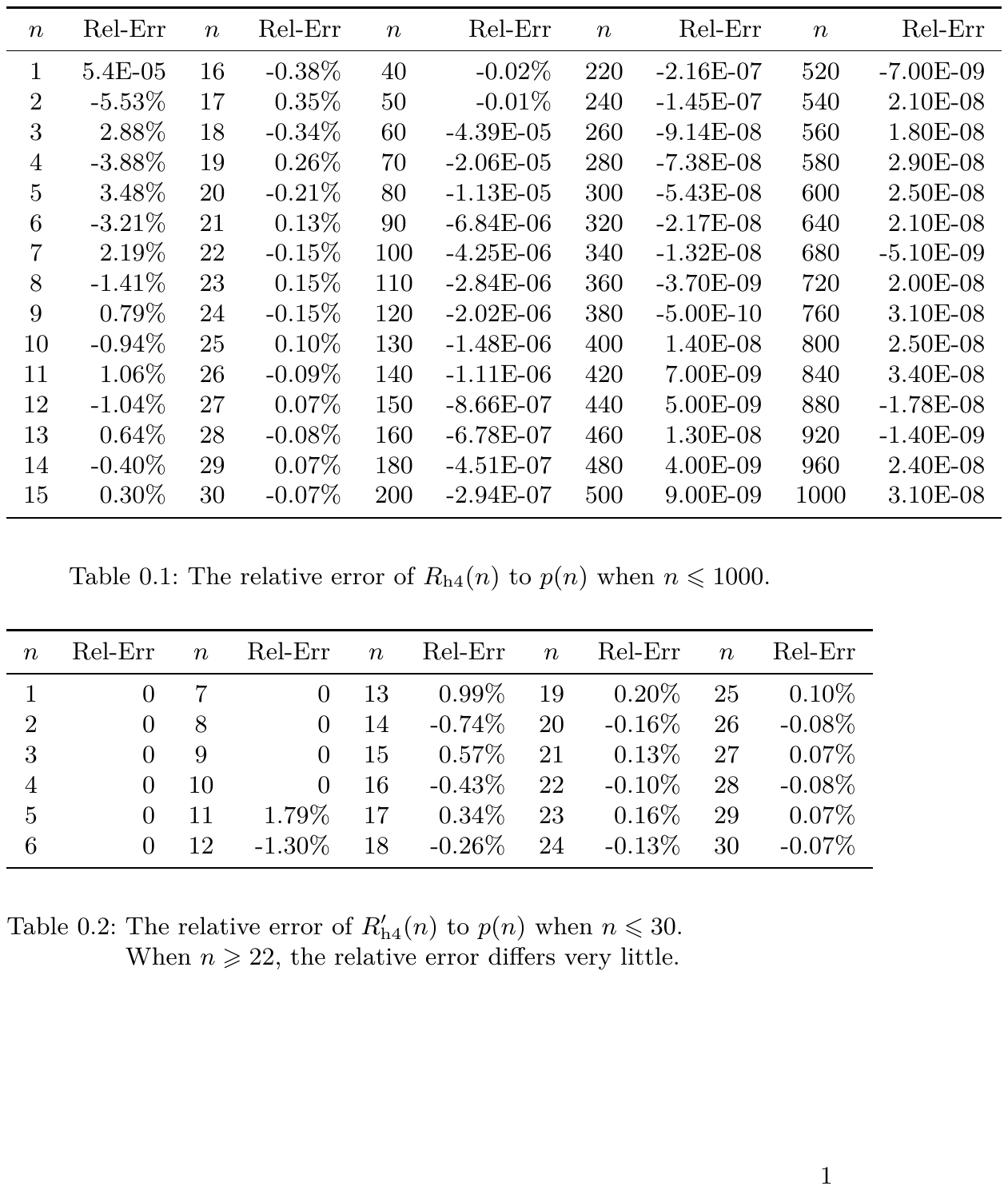}
\par\end{centering}

\caption{The relative error of $R_{\mathrm{h4}}(n)$ to $p(n)$ when $n\leqslant1000$.}
\label{Table:Rel-Err-p(n)-Rh4(n)}

\vspace{0.8cm}

\noindent \begin{centering}
\includegraphics[bb=97bp 352bp 411bp 442bp,clip,scale=0.73]{tables-v3/Table-2_4_1_f-Rel-Err-Rh4_n_-Rh4_n_round}
\par\end{centering}

\caption{The relative error of $R'_{\mathrm{h4}}(n)$ to $p(n)$ when $n\leqslant30$.\protect \\
When $n\geqslant22$, the relative error differs very little.}
\label{Table:Rel-Err-p(n)-Rh4(n)-round}

\end{table}

\section{Estimate $p(n)$ When $n\leqslant100$ \label{sec:Estmt-p(n)-less-than-100}}

Until now, all the estimation function generated for $p(n)$ are with
very good accuracy when $n$ is greater than 100, but they are not
so accurate when $n<50$. Although $R'_{\mathrm{h2}}(n)$ and $R'_{\mathrm{h4}}(n)$
are better than others, the relative error are still greater than
$1\permil$ for some values of $n$.

On the other hand, in sections \ref{sec:Fit-the-Exponent-C1(n)} and
\ref{sec:Fit-the-Denominator-C2(n)}, when $n<100$, it is nearly
impossible to fit\\
 $C_{1}\left(n\right)\doteq\dfrac{3}{2}\cdot\dfrac{\left(\ln\left(4n\sqrt{3}p(n)\right)\right)^{2}}{\pi^{2}}-n$
or\\
 $C_{2}(n)\doteq\dfrac{\exp\left(\sqrt{\frac{2}{3}}\pi\sqrt{n+C_{1}(n)}\right)}{4\sqrt{3}p(n)}-n$
\\
by a simple piecewise function with less than 4 pieces with high
accuracy, as shown on Figure \ref{Fig:2_4_3_(n_C1(n))-part-1}, Figure
\ref{Fig:2_4_4_(n_C1(n))-part-2} (on page \pageref{Fig:2_4_4_(n_C1(n))-part-2})
and Figure \ref{Fig:2_4_10-(n,C2)} (on page \pageref{Fig:2_4_10-(n,C2)}),
since the points do not lie on less than 4 smooth simple curves. 

Can we reach a better accuracy when estimating $p(n)$ by a formula
not too complicated?

In  subsection \ref{sub:Modify the Denominator only}, we fit the
data $\left(n,\ \frac{\exp\left(\pi\sqrt{\frac{2}{3}n}\right)}{4\sqrt{3}p(n)}-n\right)$
($n=20k+100$, $k$ = 1, 2, $\cdots$, 395) by a function $C_{2}(n)$
and obtained a very good estimation of $p(n)$ when $n$ > 50.

\begin{figure}[h]
\noindent \begin{centering}
\includegraphics[scale=0.26]{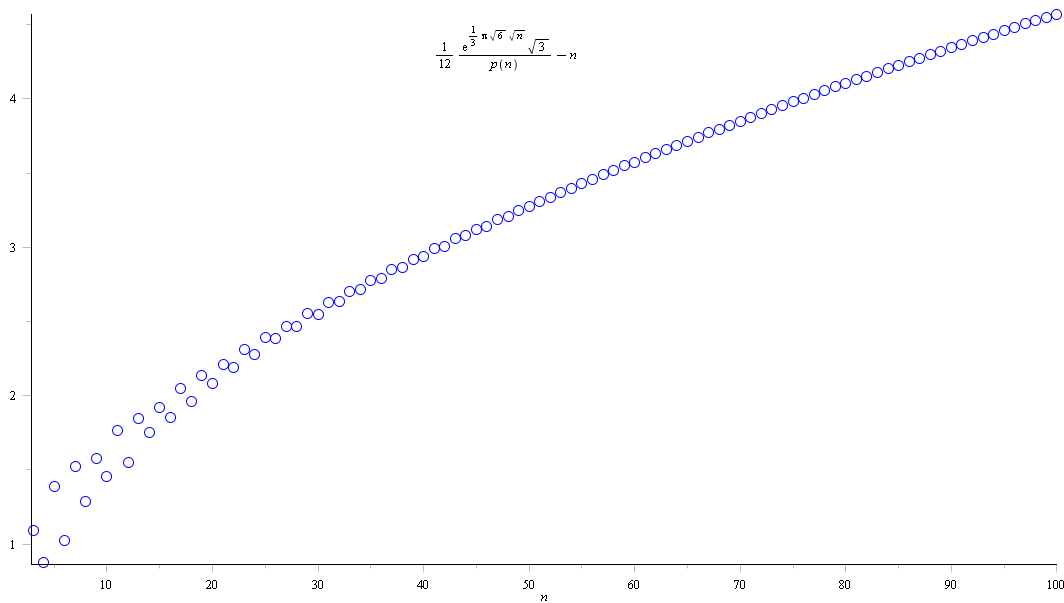}
\par\end{centering}

\caption{The graph of the data $\left(n,\,\frac{\exp\left(\sqrt{\frac{2}{3}}\pi\sqrt{n+C_{1}(n)}\right)}{4\sqrt{3}p(n)}-n\right)$
when $n$ $\leqslant$ 100}
 \label{Fig:2_4_23_(n,C2)_1-100}

\end{figure}

So we wander whether we can fit the data $\left(n,\ \frac{\exp\left(\pi\sqrt{\frac{2}{3}n}\right)}{4\sqrt{3}p(n)}-n\right)$
($n$ = 3, 4, $\cdots$, 100) by a piecewise function (with 2 pieces)
so as to get a better estimation of $p(n)$ when $n$ $\leqslant$
100?

The figure of the points of the data $\left(n,\ \frac{\exp\left(\pi\sqrt{\frac{2}{3}n}\right)}{4\sqrt{3}p(n)}-n\right)$
($n$ = 3, 4, $\cdots$, 100) are shown on Figure \ref{Fig:2_4_23_(n,C2)_1-100}
(on page \pageref{Fig:2_4_23_(n,C2)_1-100}). It is not difficult
to find that the even points (where $n$ is even) lie roughly on a
smooth curve, so are the odd points. If we try to fit them respectively,
we will have the fitting function below: 
\begin{equation}
C'_{2}(n)=\begin{cases}
0.4527092482\times\sqrt{n+4.35278}-\\
\quad\quad0.05498719946,\\
\quad\quad\quad\quad\quad\quad\quad n=3,5,7,\cdots,99;\\
0.4412187317\times\sqrt{n-2.01699}+\\
\quad\quad0.2102618735,\\
\quad\quad\quad\quad\quad\quad\quad n=4,6,8\cdots,100.
\end{cases}\label{eq:C2'(n)}
\end{equation}
 Hence we can calculate $p(n)$ by 
\begin{equation}
R_{\mathrm{h0}}(n)=\dfrac{\exp\left(\sqrt{\frac{2}{3}}\pi\sqrt{n}\right)}{4\sqrt{3}\left(n+C'_{2}(n)\right)},\quad1\leqslant n\leqslant100.\label{eq:Ram-Hardy-Rev-Fml-0A}
\end{equation}
 Consider that $p(n)$ is an integer, we can take the round approximation
of  \eqref{eq:Ram-Hardy-Rev-Fml-0A}, \label{Sym:R'h0(n)} \nomenclature[Rh0(n)]{$R'_{\mathrm{h0}}(n)$}{The Hardy-Ramanujan's revised estimation formula when $1 \leqslant n \leqslant 100$. \pageref{Sym:R'h0(n)}}
\begin{equation}
R'_{\mathrm{h0}}(n)=\left\lfloor \dfrac{\exp\left(\sqrt{\frac{2}{3}}\pi\sqrt{n}\right)}{4\sqrt{3}\left(n+C'_{2}(n)\right)}+\dfrac{1}{2}\right\rfloor ,\quad1\leqslant n\leqslant100.\label{eq:Ram-Hardy-Rev-Fml-5}
\end{equation}

The relative error of $R_{\mathrm{h0}}(n)$ (or $R'_{\mathrm{h0}}(n)$)
to $p(n)$ are shown on Table \ref{Table:Rel-Err-p(n)-Rh0(n)} (or
Table \ref{Table:Rel-Err-p(n)-Rh0(n)-round}) on page \pageref{Table:Rel-Err-p(n)-Rh0(n)}.
Compared with Table \ref{Table:Rel-Err-p(n)-HR2N(n)-round} on page
\pageref{Table:Rel-Err-p(n)-HR2N(n)-round}, we will find that when
$n\geqslant80$, $R'_{\mathrm{h2}}(n)$ is more accurate than $R'_{\mathrm{h0}}(n)$;
when $n\leqslant50$, $R'_{\mathrm{h0}}(n)$ is obviously better.

\begin{table}[h]
\noindent 

\noindent \begin{centering}
\includegraphics[bb=97bp 382bp 462bp 624bp,clip,scale=0.63]{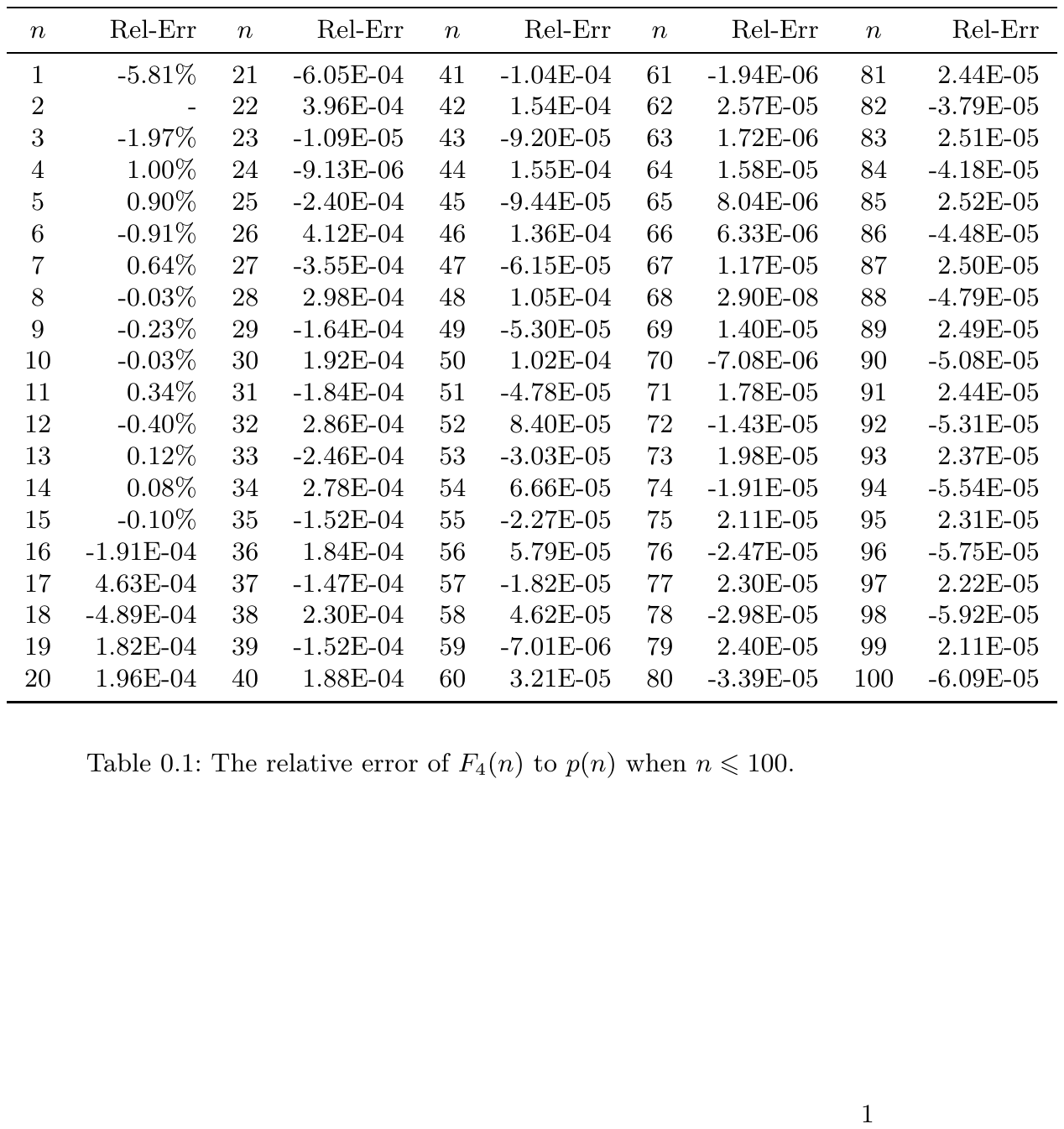}
\par\end{centering}

\caption{The relative error of $R_{\mathrm{h0}}(n)$ to $p(n)$ when $n\leqslant100$.}
\label{Table:Rel-Err-p(n)-Rh0(n)}

\vspace{0.8cm}

\noindent \begin{centering}
\includegraphics[bb=95bp 396bp 522bp 533bp,clip,scale=0.54]{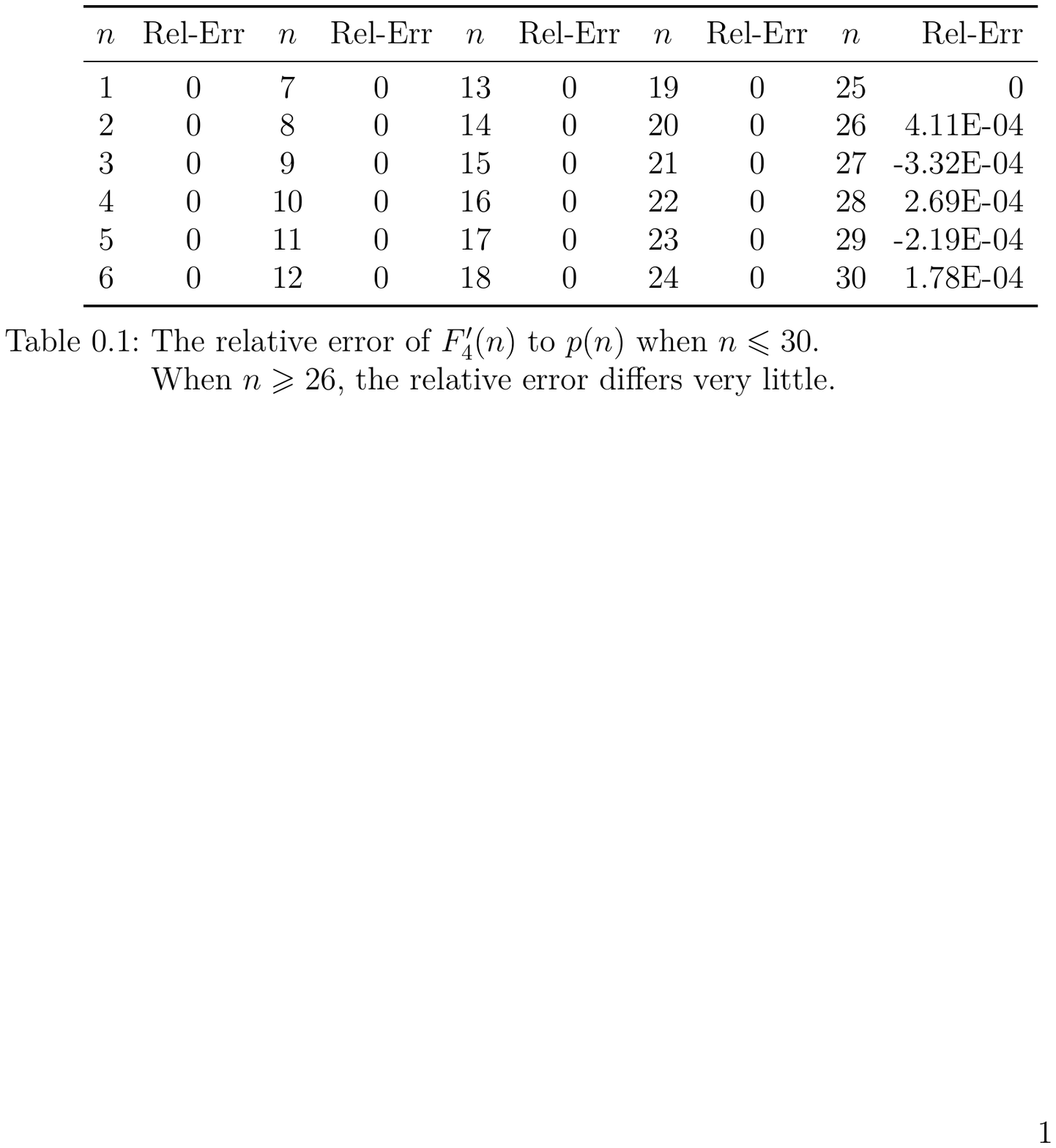}
\par\end{centering}

\caption{The relative error of $R'_{\mathrm{h0}}(n)$ to $p(n)$ when $n\leqslant30$.\protect \\
When $n\geqslant26$, the relative error differs very little.}
\label{Table:Rel-Err-p(n)-Rh0(n)-round}

\end{table}

\section{Conclusions}

In this paper, we have presented several  elementary estimation formulae
with high accuracy to calculated $p(n)$, that can be operated on
a pocket science calculator without programming function.

When $n\leqslant80$, we can use $R'_{\mathrm{h0}}(n)$ (Equation
\eqref{eq:Ram-Hardy-Rev-Fml-5}) , with a relative error less than
0.004\%; when $n>80$, we can use $R'_{\mathrm{h2}}(n)$ (Equation
\eqref{eq:Ram-Hardy-Rev-Fml-2}). 

Equations  \eqref{eq:Ram-Hardy-Rev-Fml-1}, \eqref{eq:Ram-Hardy-Rev-Fml-3}
and \eqref{eq:Ram-Hardy-Rev-Fml-4} are also very accurate although
they are not as good as  \eqref{eq:Ram-Hardy-Rev-Fml-2}.

By the construction of these estimation formulae, when $n$ $\rightarrow$
$\infty$, the relative error will approaches 0. (But the absolute
error may approaches infinity). 

If we can find the accurate expression \footnote{$\ $ such as in the form $\pi^{a/b}$ ($a,b\in\mathbb{Z}$, $ab\neq0$)
or $\pi^{a}\mathrm{e}^{b}$ ($a,b\in\mathbb{Z}$). } of the coefficients $a_{2}\doteq1.039888529$ in \eqref{eq:C4(n)},
$t_{0}\doteq0.3594143172$ in \eqref{sub:Result-2} and $a_{3}\doteq2.893270736$
in \eqref{eq:C5(n)}, and can find the explanation in theory,  
we may gain better results.

The ideas described here could be used to acquire elementary estimation
formulae in some other cases when  approximate values are frequently
wanted while the asymptotic formulae are less accurate than expectation
and the methods to calculate the exact values are inconvenient, 
such as the computation of some kinds of restricted partition numbers
if we have ( or can deduce) the asymptotic formulae beforehand.

These methods to fitting $C_{1}(n)$ and $C_{2}(n)$ could also be
used in searching for the fitting functions of some classes of data
obtain in experiments if we want more accuracy.

\section*{Acknowledgements}

\addcontentsline{toc}{section}{Acknowledgements}

The author would like to express the gratitude to his supervisor Prof.
\emph{LI Shangzhi} from Beihang University (in China) for his help
and encouragement. Special thanks must be given to Prof. \emph{Ian
M. Wanless} from Monash University (in Australia), Prof. \emph{Jack
H. Koolen} and Prof. \emph{WANG Xinmao} from USTC (in China), Dr.
\emph{ZHANG Zhe} from Xidian University and Dr. \emph{WANG Qihan}
from Anhui University (in China) for their valuable advice.

\bibliographystyle{abbrv}
\phantomsection\addcontentsline{toc}{section}{\refname}\bibliography{Ref-org-charct-22_two_Col}

$\ $
\end{document}